\newcommand{\R}{\mathds R}
\newcommand{\Z}{\mathds Z}
\newcommand{\T}{\mathds T}
\newcommand{\SL}{\mathrm{SL}}
\newcommand{\Nor}{\mathrm{Nor}}
\newcommand{\Centr}{\mathrm{Centr}}
\newcommand{\Iso}{\mathrm{Iso}}
\newcommand{\Ker}{\mathrm{Ker}}
\newcommand{\GL}{\mathrm{GL}}
\newcommand{\cF}{\mathcal F}
\newcommand{\cT}{\mathcal T}
\newcommand{\cN}{\mathcal N}
\documentclass[twoside,draft,a4paper,11pt]{amsart}

\usepackage{times}
\usepackage{dsfont}
\usepackage[all]{xy}

\usepackage{amssymb}
\numberwithin{equation}{section}

\title[Compact Lorentz manifolds with large isometry group]{On the isometry group
and the geometric structure of compact stationary Lorentzian manifolds}
\author[P.\ Piccione]{Paolo Piccione}
\address{Departamento de Matem\'atica,\hfill\break\indent
Universidade de S\~ao Paulo, \hfill\break\indent Rua do Mat\~ao
1010,\hfill\break\indent CEP 05508-900, S\~ao Paulo, SP, Brazil}
\email{piccione.p@gmail.com}
\curraddr{Departamento de Matem\'{a}ticas, \hfill\break\indent
Universidad de Murcia, Campus de Espinardo\hfill\break\indent
30100 Espinardo, Murcia, \hfill\break\indent Spain}
\author[A. Zeghib]{Abdelghani Zeghib}
\address{%
Unit\'e de Math\'ematiques Pures et Appliqu\'ees\hfill\break\indent
\'Ecole Normale Sup\'erieure de Lyon\hfill\break\indent
46, all\'ee d'Italie
\hfill\break\indent
69364 LYON Cedex 07, FRANCE}
\email{zeghib@umpa.ens-lyon.fr}
\thanks{%
The first author is partially sponsored by  Fundaci\'{o}n
S\'{e}neca grant 09708/IV2/08, Spain. Part of this work was carried out during the visit of the first author
at the \'Ecole Normale Sup\'erieure de Lyon in November 2009. The authors gratefully acknowledge the financial support
for this visit provided by \emph{ANR Geodycos}.}

\date{February 3rd, 2010}

\begin{document}

\theoremstyle{plain}\newtheorem{prop}{Proposition}[section]
\theoremstyle{plain}\newtheorem{theo}[prop]{Theorem}
\theoremstyle{definition}\newtheorem{example}{Example}
\theoremstyle{plain}\newtheorem{lem}[prop]{Lemma}
\theoremstyle{plain}\newtheorem{cor}[prop]{Corollary}
\theoremstyle{definition}\newtheorem{defin}[prop]{Definition}
\theoremstyle{remark}\newtheorem{rem}[prop]{Remark}
\theoremstyle{plain} \newtheorem{assum}[prop]{Assumption}

\theoremstyle{definition}\newtheorem*{defin*}{Definition}
\theoremstyle{plain} \newtheorem*{acknowledgement}{Acknowledgements}
\theoremstyle{definition}\newtheorem*{notation}{Notation}
\theoremstyle{remark}\newtheorem*{exa}{Example}
\theoremstyle{plain}\newtheorem*{lemn}{Lemma}
\theoremstyle{plain}\newtheorem*{theon}{Theorem}

\swapnumbers
\theoremstyle{plain}\newtheorem{theo-intro}{Theorem}
\theoremstyle{plain}\newtheorem{cor-intro}[theo-intro]{Corollary}


\begin{abstract}
We study the geometry of compact Lorentzian manifolds
that admit a somewhere timelike Killing vector field, and
whose isometry group has infinitely many connected components.
Up to a finite cover, such manifolds are products (or amalgamated
products) of a flat Lorentzian torus and a compact Riemannian (resp., lightlike) manifold.
\end{abstract}

\maketitle

\begin{section}{Introduction}

\subsubsection*{Paradigmatic example}

We will deal with dynamics and geometry of the following flavor.
Let $q$ be a Lorentz form on $\R^n$; this induces a (flat) Lorentz metric on the torus $\T^n = \R^n / \Z^n$.
The linear isometry group  of $\T^n$ is $\mathrm O(q, \Z) = \GL(n, \Z) \cap\mathrm O(q)$, and
its full isometry group is the semi-direct product $\mathrm O(q, \Z) \ltimes \T^n$.

The global and individual structure of $\mathrm O(q, \Z)$ involves interesting geometric,
arithmetic and dynamical  interactions.
For generic $q$, $\mathrm O(q, \Z)$ is trivial. Nonetheless, if $q$ is rational, i.e., if $q(x) =\sum a_{ij}x_i x_j$,
where $a_{ij}$ are rational numbers, then $\mathrm O(q, \Z)$ is \emph{big} in $\mathrm O(q)$; more precisely,
by Harich-Chandra-Borel theorem,  it is a lattice in $\mathrm O(q)$.
When  $q$ is not rational, many intermediate situations are possible.
It is a finite volume non co-compact lattice in the case of the standard form
$q_0= -x_1^2 + x_2^2 +\ldots +x_n^2$, but can be co-compact for other forms.
On the other hand, a given element $A \in\mathrm O(q_0, \Z)$, could have complicated dynamics.
For instance, if $A$ is hyperbolic, i.e., its spectrum is not contained in $\mathds S^1$, then
the remainder of eigenvalues are roots of unity or Salem numbers. Conversely, any Salem number
is the eigenvalue of such a  hyperbolic $A \in\mathrm O(q_0, \Z)$, for some dimension $n$.
\subsubsection*{Lorentz geometry and dynamics}
The global geometry of compact manifolds endowed with a non positive definite metric
(pseudo-Riemannian manifolds) can be quite different from the geometry of Riemannian manifolds.
For instance, compact pseudo-Riemannian manifolds may fail to be geodesically complete
or geodesically connected; moreover, the isometry group of a compact pseudo-Riemannian
manifold fails to be compact in general. The main goal of this paper is to investigate the
geometric structure of Lorentz manifolds \emph{essentially}  non Riemannian,
i.e., with non compact isometry group.

Lorentzian manifolds, i.e.,  manifolds endowed with metric tensors of index
$1$, play a special role in pseudo-Riemannian geometry, due to their relations with General Relativity.
The lack of compactness of the isometry group is due to the fact that,
unlike the Riemannian case, Lorentzian isometries need not be equicontinuous, and may generate
chaotic dynamics on the manifold. For instance, the dynamics of Lorentz isometries can
be of Anosov type, evocative of the fact that in General Relativity one can have
contractions of time and expansion in space.
A celebrated result of D'Ambra (see \cite{DAmbra}) states that the isometry group of a \emph{real analytic}
simply-connected compact Lorentzian manifold is compact. It is not known whether this results
holds in the $C^\infty$ case. In the last decade several authors have studied isometric group
actions on Lorentz manifolds. Most notably, a complete classification of (connected) Lie groups
that act locally faithfully and isometrically on compact Lorentzian manifolds has been obtained
independently by Adams and Stuck (see \cite{AdamsStuck}) and the second author (see \cite{Ghani1}).
Roughly speaking, the identity component $G_0$ of the isometry group of a compact Lorentz
manifold is the direct product of an abelian group, a compact semi-simple group, and, possibly,
a third factor which is locally isomorphic to either $\mathrm{SL}(2,\R)$ or to an oscillator group, or else
to a Heisenberg group.
The geometric structure of a compact Lorentz manifold that admits a faithful isometric
action of a group $G$ isomorphic to $\mathrm{SL}(2,\R)$ or to an oscillator group is well understood;
such manifolds can be described using right quotients $G/\Gamma$, where $\Gamma$ is a co-compact
lattice of $G$, and warped products, see Section~\ref{sec:identitycomponent} for more details.
Observe that such constructions produce Lorentz manifolds on which the $G_0$-action has some
timelike orbit. Recall that a Lorentz manifold is said to be \emph{stationary}
if it admits an everywhere timelike Killing vector field.
Our first result (Theorem~\ref{theorem1}) is that when the identity component of
the isometry group is non compact and it has some timelike orbit, then it must contain
a non trivial factor locally isomorphic to $\mathrm{SL}(2,\R)$ or to an oscillator group.

Thus, next natural question is to study the geometry of manifolds whose isometry
group is non compact for having an infinite number of connected components.

 \subsubsection*{Results} We will show in this paper that compact Lorentz manifolds with a \emph{large} isometry
 group
are essentially constructed in the same way. In order to define the appropriate notion
of the lack of compactness of the isometry group of a Lorentzian manifold, let us give the following:
\begin{defin*}
Let $\rho:\Gamma\to\mathrm{GL}(\mathcal E)$ be a representation of the group $G$ on the vector space $\mathcal E$.
Then, $\rho$ is said to be \emph{of Riemannian type} if it preserves some positive definite inner product on $\mathcal
E$.
We say that $\rho$ is \emph{of post-Riemannian type} if it preserves a positive semi-definite inner product
on $\mathcal E$ having kernel of dimension equal to $1$.
\end{defin*}
Observe that $\rho$ is of Riemannian type if and only if it is precompact, i.e.,
$\rho(\Gamma)$ is precompact in $\mathrm{GL}(\mathcal E)$.

Given a Lorentzian manifold $(M,\mathbf g)$, we will
denote by $\mathrm{Iso}(M,\mathbf g)$ its isometry group, and by $\mathrm{Iso}_0(M,\mathbf g)$
the identity connected component of $\mathrm{Iso}(M,\mathbf g)$. The Lie algebra of
$\mathrm{Iso}(M,\mathbf g)$ will be denoted by $\mathfrak{Iso}(M,\mathbf g)$.
By a large isometry group, we mean that the adjoint action of the discrete part
$\Gamma=\mathrm{Iso}(M,\mathbf g)/\mathrm{Iso}_0(M,\mathbf g)$ of $\mathrm{Iso}(M,\mathbf g)$ is
not of post-Riemannian type.  This implies, in particular, that
$\Gamma$ is not finite.
\smallskip

The main result of the paper is that compact Lorentz manifolds with large isometry group
are essentially built up by tori. More precisely, we prove the following structure theorem.

\begin{theo-intro}\label{thm:theorem1} Let $(M,\mathbf g)$ be a compact Lorentz manifold, and assume that
the conjugacy  action of $\Gamma$ on $\mathrm{Iso}_0(M,\mathbf g)$ is
not post-Riemannian. Then,  $\mathrm{Iso}_0(M,\mathbf g)$ contains a torus $\T= \T^d$, 
endowed with a Lorentz form $q$, such that $\Gamma $ is a subgroup of $\mathrm O(q, \Z)$.

Up to a finite cover, there is a new Lorentz metric $\mathbf g^\mathrm{new}$ on $M$  having  a larger
isometry group than $\mathbf g$, such that $\Gamma =\mathrm O(q, \Z)$.
Geometrically, $M$ is metric direct product $\T \times N$,
where $N$ is a compact Riemannian manifold, or $M$ is an amalgamated metric product
$\T \times_{\mathds S^1} L$, where $L$ is a lightlike manifold with an isometric $\mathds S^1$-action. The last
possibility holds when $\Gamma$ is a parabolic subgroup of $\mathrm O(q)$.
\end{theo-intro}
A more precise description of the original metric $\mathbf g$ is given in Section~\ref{sec:hyperbolicprodstructure}
for the hyperbolic case and in Section~\ref{sec:parabolicprodstructure} for the parabolic case.
We will in fact prove Theorem~\ref{thm:theorem1} under an assumption weaker than the non post-Rie\-mann\-ian
hypothesis for the conjugacy action of $\Gamma$. The more general statement proved here is the following:
\begin{theo-intro}\label{thm:theorem2}
Assume that $\Gamma$ is infinite and that $\Iso_0(M,\mathbf g)$ has a somewhere timelike orbit.
Then the conclusion of Theorem~\ref{thm:theorem1} holds.
\end{theo-intro}
Theorem~\ref{thm:theorem1} will follow from Theorem~\ref{thm:theorem2} once we show that, under
the assumption that the conjugacy action of $\Gamma$ is not of post-Riemannian type, then
the connected component of the identity of the isometry group must have some timelike
orbits, see Subsection~\ref{sub:Gauss}.
\smallskip

A first consequence of our main result is the following:
\begin{cor-intro}\label{thm:somewhere-everywhere}
Assume that $(M,\mathbf g)$ is a compact Lorentzian manifold with infinite discrete part $\Gamma$.
If $(M,\mathbf g)$ has a \emph{somewhere} timelike Killing vector field, then $(M,\mathbf g)$ has
an \emph{everywhere} timelike Killing vector field.
\end{cor-intro}

We will also prove (Proposition~\ref{thm:propinfmanyconncomp})
that, when $(M,\mathbf g)$ has a Killing vector field which is timelike somewhere,
then the two situations (a) and (b) below are mutually exclusive:
\begin{itemize}
\item[(a)] the connected component of the identity $\Iso_0(M,\mathbf g)$ of $\Iso(M,\mathbf g)$ is non compact;
\item[(b)] $\Iso(M,\mathbf g)$ has infinitely many connected components, as in the case of the flat Lorentzian
    torus.
\end{itemize}
The point here is that, in a compact Lorentzian manifold, the flow of a Killing vector field
which is timelike somewhere generates a non trivial precompact group in the (connected component of
the identity of the) isometry group. Thus, by continuity, the Lie algebra of the isometry group of
such manifolds must contain a non empty open cone of vectors generating precompact $1$-parameter
subgroups in the isometry group. The proof of Proposition~\ref{thm:propinfmanyconncomp}
is obtained by ruling out the existence of a non compact abelian or nilpotent factor
in the connected component of the isometry group.
The argument is based on an algebraic precompactness criterion for $1$-parameter subgroups of Lie groups
proved in Proposition~\ref{thm:charprecompactsubgr}.

Moreover, using Theorem~\ref{thm:theorem2} and previous classification results
by the second author, we prove the following partial extension of D'Ambra's result
to the $C^\infty$-realm:
\begin{theo-intro}\label{thm:simplyconnected}
The isometry group of a simply connected compact Lorentzian manifold that admits a
Killing vector field which is somewhere timelike is compact.
\end{theo-intro}

Let us give a sketch of the proof of our main result. After the preliminaries, in Section~\ref{sec:toralfactor}
we show the existence of an appropriate \emph{reduction} of the isometry group. More precisely,
we show that $\mathrm{Iso}(M,\mathbf g)$ has a closed non compact subgroup $G$ whose identity component
$G_0$ is compact abelian,  it has some timelike orbit, and the quotient $\Gamma=G/G_0$ is torsion free.
Moreover, the conjugacy action of $G$ on $G_0\cong\T^k$, denoted by
$\rho:G\to\mathrm{Out}(T^k)\cong\mathrm{GL}(k,\Z)$,
preserves some Lorentz form of $\R^k$; thus, the image $\rho(G)$ is contained in
$\mathrm{GL}(k,\Z)\cap\mathrm{SO}(q)$
for some Lorentz quadratic form $q$.
The proof of our structure result is then obtained by relating the dynamics of an isometry
$f\in G$, $f\not\in G_0$, on $M$, and the dynamics of the conjugacy action of $f$ on the toral
factor $G_0$. In this step the proof splits into two cases, according to whether $\rho(f)$ is a parabolic
or a hyperbolic element of $\mathrm{SO}(q)$. We will make full use of the techniques developed by the
second author in \cite{Zeghib.GAFA1}, in particular we will employ the notions of \emph{approximate stable foliation}
and of \emph{strongly approximate stable foliation} of an isometry of a compact Lorentz manifold.
This analysis will show the existence of a possibly smaller torus $\T^d\subset\T^k$, $d\ge2$,
whose action on $M$ is everywhere free and timelike. In the hyperbolic case, the orthogonal distribution
to the $\T^d$-orbits is integrable; a general covering Lemma (Section~\ref{sec:generalcovering})
will yield the product structure of (a finite covering of) $M$. The amalgamated product structure
in parabolic case is somewhat more involved,
and it is discussed in Section~\ref{sec:parabolicprodstructure}.
\end{section}
\begin{section}{Preliminaries}
In this section we will collect several auxiliary results needed in the rest of the paper.
\subsection{Toral subgroups}
Our aim here is to determine the freeness of isometric toral actions on manifolds.
The key fact is that the set $S(\T^d)$ of all closed subgroups of the $d$-torus $\mathds T^d$
is countable, and it satisfies a uniform discreteness property.
\begin{lem}\label{thm:prelimisotropy}
Let $X$ be a locally compact metric space, and let $\phi:X \to S(\mathds T^d)$ be a \emph{semi-continuous} map,
that is, if $x_n \to x$, then any limit of $\phi(x_n)$ is contained in $\phi(x)$. Then, there
exists $A\in S(\mathds T^d)$ such that $\phi^{-1}(A)$ has non empty interior.
\end{lem}
\begin{proof}
For $A \in S(T)$, set $F_A = \big\{x \in X:\text{such that}\ A\subset\phi(x)\big\}$. By the semi-continuity,
the closure $\overline{\phi^{-1}(A)}\subset F_A$ for all $A\in S(\mathds T^d)$. Clearly,
$X=\bigcup\limits_{A\in S(\mathds T^d)}\phi^{-1}(A)$.
By Baire's theorem, the interior $\mathrm{int}\big(\overline{\phi^{-1}(A)}\big)$ of some $\overline{\phi^{-1}(A)}$
must be non empty. Thus, the intersection $\mathrm{int}\big(\overline{\phi^{-1}(A)}\big)\cap\phi^{-1}(A)$ is
non empty. Let $x$ be a point of such intersection, so that $A=\phi(x)$, and there is a neighborhood $V$ of $x$
such that  $\phi(y)\supset \phi(x)$ for all $y\in V$. By semi-continuity, we must have
equality $\phi(y)=\phi(x)$ for $y$ in some neighborhood $V'\subset V$.
This follows from the fact that $A$ is an isolated point of the set:
\[S(\T^d;A)=\big\{B\in S(\mathds T^d):A\subset B\big\},\]
see Lemma~\ref{thm:isolatedpoint}.
Hence, $\phi^{-1}(A)$ has non empty interior.
\end{proof}
\begin{lem}\label{thm:isolatedpoint}
Every $A\in S(\T^d)$ is isolated in $S(\T^d;A)$.
\end{lem}
\begin{proof}
Let us consider the case that $A$ is the trivial subgroup. To prove that $A=\{1\}$ is isolated
in $S(\T^d)$ it suffices to observe that there exist two disjoint closed subsets $C_1,C_2\subset\T^d$
such that $C_1$ is a neighborhood of $1$, $C_1\cap C_2=\emptyset$, and with the property that if $B\in S(\T^d)$
is such that $B\cap C_1\ne\{1\}$, then $B\cap C_2\ne\emptyset$. For instance, one can take $C_1$ to be the closed
ball
around $1$ of radius $r>0$ small, and $C_2=\big\{p\in\T^d:2r\le\mathrm{dist}(p,1)\le3r\big\}$. Here we are
considering the distance on $\T^d=\R^d/\Z^d$ induced by the Euclidean metric of $\R^d$. The
number $r$ can be chosen in such a way that $C_1\cap C_2=\emptyset$. Every non trivial element
in $C_1$ has some power in $C_2$, which proves that $C_1$ and $C_2$ have the required properties.
Thus, if $A_n\in S(\T^d)$ is any sequence which is not eventually equal to $A$,
then some limit point of $A_n$ must be contained in $C_2$, and therefore $\lim A_n\ne A$.
If one replaces $\T^d$ by a finite quotient of $\T^d$, then one gets to the same conclusion by essentially
the same proof. The general case is obtained by considering the quotient $\T^d/A$, which is equal to a finite
quotient of a torus.
\end{proof}
\begin{cor}\label{thm:countableclosedsubgroup}
If $X$ is a locally compact metric space and $\phi:X\to S(\mathds T^d)$ is semi-continuous,
then, there is a dense open subset $U \subset X$, where $\phi $ is locally constant, i.e., any
$x \in U$ has a neighborhood $V$ where $\phi$ is constant.
\end{cor}
\begin{proof}
Let $U$ be the open subset of $X$ given by the union of the interiors of the sets $\phi^{-1}(A)$,
with $A$ running in $S(\mathds T^d)$. This is the largest open subset of $X$ where $\phi$ is locally
constant. If $U$ were not dense, then there would exist a non empty open subset $V\subset X$
with $V\cap U=\emptyset$. The restriction $\widetilde\phi$ of $\phi$ to $V$ is a semi-continuous map, with the
property
that $\widetilde\phi^{-1}(A)$ has empty interior for all $A\in S(\mathds T^d)$. By Lemma~\ref{thm:prelimisotropy},
this is impossible, hence $U$ is dense.
\end{proof}
\begin{cor}\label{thm:faithfulfreeaction}
Any faithful isometric action of a torus $\T^d$ on some pseudo-Rie\-mann\-ian manifold $(M,\mathbf g)$
is free on a dense open subset of $M$.
\end{cor}
\begin{proof}
Apply Corollary~\ref{thm:countableclosedsubgroup} to the map $\phi:M\to S(\mathds T^d)$ that
associates to each $p\in M$ its stabilizer. Such map is obviously semi-continuous. Thus,
on a dense open subset $U$ of $M$, the stabilizer of the isometric action is locally constant.
No nontrivial isometry of a pseudo-Riemannian manifold fixes all points of a non empty open subset,
and this implies that the stabilizer of each point of $U$ is trivial.
\end{proof}

\subsection{Linear dynamics}

Gauss maps (and variants) have the advantage to transform the dynamics on $M$ into a linear dynamics, i.e., an action
of the group in question on a linear space or an associated projective space, via a linear representation. We will
prove in the sequel a stability result: if a linear group ``almost-preserves'' a Lorentz form, then it (fully)
preserves  another one.
We start with the individual case, i.e., with  actions of the infinite cyclic group $\Z$, and then we will
consider general groups.

\subsubsection{Individual dynamics}
Let $\mathcal E$ be a vector space, and $A \in \mathrm{GL}(\mathcal E)$. It has a Jordan decomposition $A =EHU$,
where
$U$ is unipotent (i.e., $U-1$ is nilpotent), $H$ hyperbolic (i.e., diagonalizable over $\R$), and $E$ is elliptic
(i.e.,
diagonalizable over $\mathds C$, and all its eigenvalues have norm equal to $1$).

If $F$ is a space obtained from $\mathcal E$ by functorial constructions, e.g., $F =\mathrm{Sym}(\mathcal E^*)$
the space of quadratic forms on $\mathcal E$, or $F=\mathrm{Gr}^d(\mathcal E)$  the Grassmannian of $d$-di\-mensional
subspaces of $\mathcal E$, the associated $A$-action on $F$  will be denoted by $A^F$.
Naturally, when $F$ is a vector space, we have $A^F = E^F H^F U^F$.

A point $p \in \mathcal E$ is \emph{$A$-recurrent} if there is
$n_i\in\mathds N $, $n_i \to \infty$, such that $A^{n_i}(p) \to p$ as $i\to\infty$.
A point $p$ is \emph{$A$-escaping} if for any compact subset $K\subset \mathcal E$ there is $N$ such that
$A^n(p)\notin K$, for $n>N$. So, $p$ is non-escaping if there is $n_i \to \infty$,
such that $A^{n_i} p$ stay in some compact set $K \subset \mathcal E$.

It is easy to prove the following
\begin{lem}\label{thm:prelimlineardynamics1}
Let $p \in F$ be a point recurrent under the $A^F$-action. Then $p$ is fixed by
$H^F$ and $U^F$.
If $F$ is a vector space, $p$ is $A^F$-non-escaping iff $p$ is fixed by $H^F$ and
$U^F$.
\end{lem}

Recall that an element $A$ of the orthogonal group of a Lorentz form $q$ is  \emph{hyperbolic}
if it has one eigenvalue $\lambda$ with $\vert\lambda\vert\neq1$, and $A$   is \emph{parabolic} if it is not
diagonalizable  (over $\mathds C$).
In other words, $A$ is hyperbolic  if it is conjugate in $\GL(k,\R)$  (in fact in $\mathrm{SO}(q)$)
to a matrix of the form:
\[\begin{pmatrix}\begin{matrix}\lambda&0\\0&\lambda^{-1}\end{matrix}&\mathbf0\\\mathbf0&R\end{pmatrix}\]
where $\lambda\in\R$, $\lambda>1$, and $R\in\mathrm{SO}(k-2)$.
Similarly, $A$ is parabolic if it has the normal form:
\[\begin{pmatrix}\begin{matrix}1&t&-\tfrac{t^2}2\\0&1&t\\0&0&1\end{matrix}&\mathbf0\\\mathbf0&R\end{pmatrix}\]
with $t\in\R$ and $R\in\mathrm{SO}(k-3)$.
\begin{lem} \label{average}
Set $F =\mathrm{Sym}(\mathcal E)$, and assume $A=EHU$  non-elliptic (i.e., either $H$ or $U$
is non-trivial). Suppose there is  a Lorentz form $q_0$  which is $A$-recurrent, and
let $K\subset\mathrm{GL}(\mathcal E)$ be the torus generated by the powers of $E$. Then,
$\int_K B^F(q_0)\,\mathrm d\mu(B)$ is an $A$-invariant  Lorentz form, where $\mu$ is the Haar measure on $K$.
\end{lem}
\begin{proof}
Assume $H\ne1$. By Lemma~\ref{thm:prelimlineardynamics1}, $H$ preserves $q_0$.
Then, there exists $\lambda\in\R\setminus\{1\}$, such that $\lambda,\lambda^{-1}\in\sigma(A)$.
The eigenspaces $V_\lambda$ and $V_{\lambda^{-1}}$ are isotropic, hence $1$-dimensional. Namely,
for $v\in V_\lambda$, $q_0(v,v)=q_0(Hv,Hv)=\lambda^2q_0(v,v)$; similarly, for $v\in V_{\lambda^{-1}}$,
$q_0(v,v)=\lambda^{-2}q_0(v,v)$. In both cases, one must have $q_0(v,v)=0$.
It follows that the direct sum
$\mathcal E_\lambda=V_\lambda\oplus V_{\lambda^{-1}}\cong\R^2$ is timelike, and the
Lorentzian form $q_0\vert_{\mathcal E_\lambda}$ is $H$-invariant.

Observe now that $\mathcal E_\lambda$ is the unique 2-space on which $H$ is non-elliptic
(i.e., all its powers are uniformly unbounded). Thus, any endomorphism commuting with $H$ preserves $\mathcal
E_\lambda$.
Therefore, $E$ preserves $\mathcal E_\lambda$. A similar argument yields that $E$ preserves the orthogonal $\mathcal
E_\lambda^\perp$ (with respect to $q_0$), since it is the maximal subspace on which $H$ acts elliptically.

On the other hand, up to multiples, the only quadratic form on $\R^2$ which is preserved by $\begin{pmatrix}
\lambda&0\cr0&\lambda^{-1}\end{pmatrix}$,
with $\lambda\ne1$, is the form $(x_1,x_2)\mapsto x_1x_2$. As above, since $E$ commutes with $H$, it preserves this
quadratic form on $\mathcal E_\lambda$.

From all this, it follows that $E^n$ preserves the fixed Lorentzian metric $q_0\vert_{\mathcal E_\lambda}$
for all $n$, while in the orthogonal space to $\mathcal E_\lambda$
the restriction of $q_0$ is positive definite. Hence, the average $\int_K B^F(q_0)\,\mathrm d\mu(B)$,
which is $A$-invariant, is equal to $q_0$ on $\mathcal E_\lambda$ and positive definite on the orthogonal
of $\mathcal E_\lambda$.

Now the case $H= 1$ but $U\ne1$ can be treated in analogous way, using essentially that $E$ commutes with $U$, and
that $U$ has the normal form of a parabolic element in $\mathrm{SO}(q_0)$, see the proof of Lemma
\ref{hyperbolic.parabolic} for a similar argument.
\end{proof}

\subsubsection{Group dynamics}
We consider now a group $\Gamma$ acting on $\mathcal E$ via a representation $\rho: \Gamma \to
\GL(\mathcal E)$. One can naturally develop a theory of recurrence  leading in particular to  a
variant  of the previous lemma for groups.
\begin{prop} \label{quasi.preserve} Let $ \rho: \Gamma \to \GL(\mathcal E)$ be  such that $\rho(a)$ is non-elliptic
for
any $a \in \Gamma$.  Let $F =\mathrm{Sym}(\mathcal E)$, and assume that the associated action
$\rho^F$ preserves a compact set of $F$ contained in the (open) subset of Lorentz forms, and that $\rho^F$ leaves
invariant a finite measure on such compact set.  Then,   $\rho(\Gamma)$  preserves some Lorentz form.
\end{prop}

The proof of Proposition~\ref{quasi.preserve} occupies the remainder of this subsection.
\smallskip

For $a \in \Gamma$, let $H(a)$ and $U(a)$ be the hyperbolic and
unipotent parts of $\rho(a)$, respectively.
Let $q_0$ be a form recurrent under the $\Gamma$-action. By Lemma \ref{thm:prelimlineardynamics1},
$q_0$ is preserved by any $H(a)$ and any $U(a)$, $a \in \Gamma$.
Let $\Gamma^\star$ be the group generated by $\{H(a),\ U(a): a \in \Gamma \}$. Then $\Gamma^\star$ is a subgroup
of $\mathrm{SO}(q_0)$.

A subgroup  $L \subset\mathrm{SO}(q_0)$ is  non-elementary
 if it  does not preserve a timelike or a lightlike direction  in $\mathcal E$.
To simplify, let us assume that $\Gamma^\star$ is not elementary;
the elementary case can be treated separately.
For instance, a subgroup of $\mathrm{SO}(q)$ preserves a timelike direction if and only if
it is pre-compact. More precisely, preserving a Lorentz form and a timelike direction
is equivalent to preserving a positive definite inner product.\footnote{
If a subgroup $H$ of $\mathrm{SO}(q)$ preserves a direction spanned by a timelike vector $e$,
then it also preserves its spacelike orthogonal, and therefore it preserves the positive product
obtained from $q$ by reversing the sign of $q(e)$. Conversely, if $H$ preserves both a Lorentz form
$q$ and a positive definite inner product $q_0$, then it commutes with the $q_0$-symmetric operator
$S$ defined by $q=q_0(S\cdot,\cdot)$. Such operator $S$ has precisely one negative eigenvalue, and
the corresponding timelike direction is preserved by $H$. Note that this implies in particular that
no precompact subgroup of $\mathrm{SO}(q)$ acts irreducibly on $\mathcal E$.}
The case of a lightlike direction looks like the  previous case of a group generated by a single element.
\begin{lem}\label{thm:decomposition} Let $L$ be a non-elementary subgroup of $\mathrm{SO}(q_0)$. Then:
\begin{enumerate}
\item\label{itm:lindyn1} There is a unique decomposition $\mathcal E = \mathcal A \oplus \mathcal B$,
characterized by the fact that $H$ acts precompactly on $\mathcal B$ and
non precompactly on  $\mathcal A$.
Furthermore, the decomposition is $q_0$-orthogonal, with  $\mathcal A$  timelike and
 $\mathcal B$ spacelike.
\item\label{itm:lindyn2} Any  Lorentz form $q$ invariant under $L$ decomposes as $\lambda {q_0} \vert_{{\mathcal
    A} } + \delta$,
where $\delta$ is a positive definite form on $\mathcal B$.
\end{enumerate}
It follows in particular that any $A \in\mathrm{GL}(\mathcal E)$ normalizing $L$ respects the decomposition, and the
image $A^*q_0$  has the previous form.
\end{lem}
\begin{proof}
For part \eqref{itm:lindyn1}, if $L$ does not act irreducibly, then it
leaves invariant  some subspace $R \subset \mathcal E$. By
the non-elementarity hypothesis, $R$ is spacelike, or timelike with dimension larger than $1$.
 If $R$ is spacelike, then we consider its orthogonal $R^\perp$ which is timelike and with dimension greater than
 $1$.
 Iterating the process, we get a timelike subspace $\mathcal A$ of dimension greater than $1$, having no proper
 $L$-invariant subspace, that is, $L$ acts irreducibly on $\mathcal A$. One can show that
 the orbit $Lu$ is unbounded for any non zero
 vector $u \in \mathcal A$.
 Namely, let $\mathcal A_0$ be the subspace of $\mathcal A$ consisting of vectors $u$ whose $L$-orbit is bounded.
 This is clearly $L$-invariant, and by irreducibility it is either $0$ or the whole $\mathcal E$.
 But $\mathcal A_0$ cannot coincide with $\mathcal A$, because otherwise the image of $L$ in
 $\mathrm{SO}\big(q\vert_\mathcal A\big)$
 would be precompact, and no precompact subgroup of $\mathrm{SO}\big(q\vert_\mathcal A\big)$ acts irreducibly.
 Thus, every non zero element of $\mathcal A$ has unbounded $L$-orbit.

 Clearly, the action of $L$ on the orthogonal space $\mathcal B=\mathcal A^\perp$
 is precompact, as $q_0$ is positive definite on $\mathcal B$.

 Let $(u, v) \in \mathcal A \oplus \mathcal B$, with $u \neq 0$. Then, there is a sequence $l_i \in L$,
 such that the direction of $l_i(u, v)$ converges to a direction in $\mathcal A$, because
 $l_i(u)$ is unbounded and $l_i(v)$ is bounded.
 Hence, $\mathcal A$ is an attractor, and such dynamical characterization
 implies the uniqueness of $\mathcal A$ as an irreducible timelike $L$-invariant subspace.

For part \eqref{itm:lindyn2}, let $q$ be another $L$-invariant Lorentz form.
Using the dynamical properties of $L$ on $\mathcal A$ and $\mathcal B$, one proves that
$\mathcal A$  and $\mathcal B$ are $q$-orthogonal.
Namely, consider the set $\{u\in\mathcal A:u\ \text{is $q$-orthogonal to } \mathcal B\}$.
This is an $L$-invariant subspace of $\mathcal A$, and by irreducibility it is either $\{0\}$
or equal to $\mathcal A$. Let us show that it is non zero. Fix $u\in\mathcal A$, $u\ne0$, and
let $l_i$ be a sequence in $L$ such that $l_iu$ is unbounded. Choose an auxiliary norm $\Vert\cdot\Vert$
in $\mathcal A$ and set $u_i=l_iu/\Vert l_iu\Vert$; up to subsequences, we can assume that
$\lim_iu_i=u_\infty\in\mathcal A$
and $u_\infty\ne0$. We claim that $u_\infty$ is $q$-orthogonal to $\mathcal B$; namely, given $v\in\mathcal B$,
$q(u_\infty,v)=\lim_iq(u_i,v)=\lim_iq(u,l_i^{-1}v)/\Vert l_iu\Vert=0$, because $l_iu$ is bounded.

It remains to check the proportionality condition
along $\mathcal A$. For this, we can assume $\mathcal B=0$, and that $L$ acts irreducibly on $\mathcal E$.
Then, $q$ can be written by means of $q_0$ via an endomorphism $S$ of $\mathcal E$, i.e.,
as bilinear forms, $q=q_0(S\cdot,\cdot)$. The fact that $L$
preserves both $q$ and $q_0$ implies that $S$ commutes with the elements of $L$.
Let us show that $S$ has at least one real eigenvalue. Consider the group
$R=\mathrm{SO}(q)\cap\mathrm{SO}(q_0)$; since $R$ contains $L$, then it acts irreducibly on
$\mathcal E$. Now, $R$ must contain either a hyperbolic or a unipotent element.
Namely, $R$ is a non compact algebraic group, and it has finitely many connected components.
The connected component of the identity of $R$ is a connected non compact Lie group, so it must contain
a non precompact $1$-parameter subgroup. Hence, it contains at least one non elliptic element.
Also, an algebraic group is closed under Jordan decomposition, i.e., if $A\in R$, then its elliptic, its hyperbolic
and
its unipotent factors belong to $R$. This proves that $R$ contains either a hyperbolic or a unipotent element.
Such element has at least one real eigenvalue of multiplicity $1$, and the corresponding eigenline is preserved
by $S$. Thus, $S$ has at least one real eigenvalue. The corresponding $S$-eigenspace is $L$-invariant, and by
irreducibility, such eigenspace must coincide with $\mathcal E$.
It follows that $S$ is a homothety, and this concludes the proof.
\end{proof}
Apply Lemma~\ref{thm:decomposition} to $L = \Gamma^\star$,    with the key observation that
$\Gamma$ normalizes $\Gamma^\star$,
that is,
$\rho(a) \Gamma^H   \rho(a)^{-1} =  \Gamma^H  $ for all $a$, since $H(aba^{-1}) = \rho(a) H(b) \rho(a)^{-1}$ and
$U(aba^{-1}) = \rho(a) U(b) \rho(a)^{-1}$. It follows that $\Gamma$ preserves the decomposition,
and $\rho(a)^*q_0 = \lambda (a) {q_0}\vert_{\mathcal A} + \delta (a)$. We know that
$H(a)$ and $U(a)$ preserve $q_0$, and thus the image by $\rho(a)$ coincides with its  image under  the elliptic part
$E(a)$,
that is, $E(a)$ acts homothetically on $(\mathcal A, q_0)$. But, since $E(a)$ is elliptic,  we have, $\lambda (a) =
\pm 1$.
In fact, $\lambda(a) = +1$, since we know that $\rho^*(a) q_0$ is a Lorentz form.

Summarizing, the action of $\Gamma$ is $q_0$-isometric on $\mathcal A$, and each element
  $a$ acts on $\mathcal B$ via its elliptic part $E(a)$, in other words ${\rho(a)}\vert_{\mathcal B} =
  {E(a)}\vert_{\mathcal B}$. In order to prove that $\Gamma$ preserves a Lorentz form
  on $\mathcal E$ it suffices to show that its action (via the $E(a)$) preserves a positive definite form on
  $\mathcal B$.  More  precisely,
  let $\Gamma^E$ be the subgroup of $\mathrm{GL}(\mathcal B)$ generated by
  the $E(a)$'s. It is true that any $A \in \Gamma^E$ is elliptic, in particular it preserves a positive definite
  form.
  However, there are exotic examples of groups with all elliptic elements that are not
  pre-compact (see \cite{Waterman})! This is why in the statement of Proposition~\ref{quasi.preserve} we make an
  assumption
  stronger than just recurrence, but rather the existence of a an invariant measure.
  The hypothesis can be here seen as equivalent to the existence of an invariant measure on the space
  $\mathrm{Sym}(\mathcal B)$ supported on the open subset of positive-definite quadratic forms.
  At this stage, one could use Furstenberg Lemma which says, roughly speaking, that linear groups preserving
  a volume act pre-compactly on its support, see \cite{Zim}.
 The proof of this statement in the present situation is straightforward. If $\mu$ is the preserved measure, then the
 mean
  $\int_{\mathrm{Sym}(\mathcal B)} q\, \mathrm d\mu(q)$ is a positive definite quadratic form, invariant under the
  $\Gamma$-action.

This completes the proof of Proposition \ref{quasi.preserve}. $\Box$

\begin{cor}\label{thm:finiteindexsubgroup}
Let $\Gamma$ be a subgroup of $\GL(k,\Z)$ which acts on
$\mathrm{Sym}(\R^k)$ by preserving a finite measure supported in a
  in the open set of Lorentz forms. Then, up to a finite index, $\Gamma$ preserves a Lorentz form.
\end{cor}
\begin{proof} Indeed, up to a finite index, we can assume that $\Gamma$ has no torsion elements.
This follows from Selberg Lemma (which says that a finitely generated matrix group has  a torsion free subgroup of
finite index, see for instance \cite{Roger}). But, in  a discrete group an elliptic $A$ element has finite order,
because the group $\{ A^n, n \in  \Z\}$ is bounded and thus finite. We deduce that, up to a finite index, $\Gamma$
has no elliptic element, and thus Proposition~\ref{quasi.preserve} applies.
\end{proof}

\subsection{A Gauss map}\label{sub:Gauss}
Let $(M,\mathbf g)$ be a compact Lorentzian manifold,  let $\Iso(M,\mathbf g)$ denote its isometry group, which is Lie
group (see
for instance \cite{Koba}), and denote by $\Iso_0(M,\mathbf g)$ the connected component
of the identity of $\Iso(M,\mathbf g)$. The Lie algebra
of $\Iso(M,\mathbf g)$ will be denoted by $\mathfrak{Iso}(M,\mathbf g)$; let us recall that there is a Lie algebra
anti-isomorphism from $\mathfrak{Iso}(M,\mathbf g)$ to the space of Killing vector field $\mathrm{Kill}(M,\mathbf g)$
obtained by mapping a vector $\mathfrak v\in\mathfrak{Iso}(M,\mathbf g)$ to the Killing field $K^{\mathfrak v}$
which is the infinitesimal generator of the one-parameter group of isometries $\R\ni t\mapsto\exp(t\mathfrak
v)\in\Iso(M,\mathbf g)$.
If $\Phi$ is a diffeomorphism of $M$ and $K$ is a vector field on $M$, we will denote by $\Phi_*(K)$ the
\emph{push-forward} of $K$ by $\Phi$, which is the vector field given by
$\Phi_*(K)(p)=\mathrm d\Phi\big(\Phi^{-1}(p)\big)K\big(\Phi^{-1}(p)\big)$ for all $p\in M$.
If $\Phi$ is an isometry and $K$ is Killing, then $\Phi_*(K)$ is Killing.

If $\Phi\in\Iso(M,\mathbf g)$, then:
\begin{equation}\label{eq:equivariance}
\phantom{\quad\forall\,\mathfrak v\in\mathfrak{Iso}(M,\mathbf g).}\Phi_*(K^{\mathfrak
v})=K^{\mathrm{Ad}_\Phi(\mathfrak v)},\quad\forall\,\mathfrak v\in\mathfrak{Iso}(M,\mathbf g).
\end{equation}
It will be useful to introduce the following map.
Let $\mathrm{Sym}(\mathfrak g)$ denote the vector space of symmetric bilinear forms on $\mathfrak g$.
The \emph{Gauss map} $\mathcal G:M\to\mathrm{Sym}(\mathfrak g)$ is the map defined by
\begin{equation}\label{eq:Gaussmap}
\mathcal G_p(\mathfrak v,\mathfrak w)=\mathbf g_p\big(K^\mathfrak v(p),K^\mathfrak w(p)\big),
\end{equation}
for $p\in M$ and $\mathfrak v,\mathfrak w\in\mathfrak g$. The following identity is immediate:
\begin{equation}\label{eq:equivGauss}
\mathcal G_{\Phi(p)}=\mathcal G_p\big(\mathrm{Ad}_\Phi\cdot,\mathrm{Ad}_\Phi\cdot\big),
\end{equation}
for all $\Phi\in\mathrm{Iso}(M,\mathbf g)$. In this paper we will be interested in the case
where $(M,g)$ admits a Killing vector field which is timelike somewhere. In this situation,
the image of the Gauss map contains a Lorentzian (nondegenerate) symmetric bilinear form
on $\mathfrak g$ (in fact, a non empty open subset consisting of Lorentzian forms,
this will be used in Lemma~\ref{thm:preserveLorntzform}).

We now have the necessary ingredients to show how the proof of Theorem~\ref{thm:theorem1}
is obtained from Theorem~\ref{thm:theorem2}.
\begin{proof}[Proof of Theorem~\ref{thm:theorem1} from Theorem~\ref{thm:theorem2}]
Let us assume that the action of the discrete part $\Gamma$ on $\mathrm{Iso}_0(M,\mathbf g)$ is not of post-Riemannian
type;
we will show by contradiction that  $\mathrm{Iso}_0(M,\mathbf g)$ has a somewhere timelike orbit.
Let $\kappa$ be the quadratic form on $\mathfrak{Iso}(M,\mathbf g)$ defined by
\[\kappa(\mathfrak v,\mathfrak w)=\int_M\mathcal G_p(\mathfrak v,\mathfrak w)\,\mathrm dM(p),\]
the integral being taken relatively to the volume element of the Lorentzian metric
$\mathbf g$. By \eqref{eq:equivariance}, $\kappa$ is invariant by the conjugacy action.
If $\mathrm{Iso}_0(M,\mathbf g)$ has no timelike orbit, then $\kappa$ is positive semi-definite.
The proof will be concluded if we show that the kernel $\mathcal K=\mathrm{Ker}(\kappa)$ has dimension less than or
equal to $1$.  Assume that $\mathcal K$ is not trivial, i.e., that $\kappa$ is positive
semi-definite. If $\mathfrak v\in\mathcal K$, then for all $\mathfrak w\in\mathfrak{Iso}(M,\mathbf g)$
and all $p\in M$, $\mathbf g_p(K^\mathfrak v_p,K^\mathfrak w_p)=0$; in particular,
$\mathbf g_p(K^\mathfrak v_p,K^\mathfrak v_p)=0$, i.e., $K^\mathfrak v$ is an everywhere isotropic\footnote{%
Here we use the following terminology: a vector $v\in TM$ is \emph{isotropic} if $\mathbf g(v,v)=0$, and
it is \emph{lightlike} if it is isotropic and non zero.}
Killing vector field on $M$. Non trivial isotropic Killing vector fields are never vanishing,
see for instance \cite[Lemma~3.2]{BeeEhrMar}; this implies that the map $\mathcal K\ni\mathfrak v\mapsto K^\mathfrak
v_p\in T_pM$
is an injective vector space homomorphism for all $p\in M$.
On the other hand, its image has dimension $1$, because an isotropic subspace of a Lorentz form
has dimension at most $1$, hence $\mathcal K$ has dimension $1$.
\end{proof}
\subsection{Precompactness of $1$-parameter sugroups}
It will be useful to recall that there is a natural smooth left action of $\Iso(M,\mathbf g)$ on the principal bundle
$\mathcal F(M)$ of all linear frames of $TM$, defined as follows.
If $b=(v_1,\ldots,v_n)$ is a linear basis of $T_pM$, then for $\Phi\in\Iso(M,\mathbf g)$ set
$\overline\Phi(b)=\big(\mathrm d\Phi_p(v_1),\ldots,\mathrm d\Phi_p(v_n)\big)$, which is a basis of $T_{\Phi(p)}M$.
The action of $\Iso(M,\mathbf g)$ on $\mathcal F(M)$ is defined by $\Iso(M,\mathbf g)\times\mathcal
F(M)\ni(\Phi,b)\longmapsto\overline\Phi(b)\in\mathcal F(M)$.
Given any frame $b\in\mathcal F(M)$,
then the map $\Iso(M,\mathbf g)\ni\Phi\mapsto\overline\Phi(b)\in\mathcal F(M)$ is an embedding of $\Iso(M,\mathbf g)$
onto a closed
submanifold of $\mathcal F(M)$ (see \cite[Theorem~1.2, Theorem~1.3]{Koba}), and thus the topology and the
differentiable structure of $\Iso(M,\mathbf g)$ can be studied by looking at one of its orbits in the frame bundle.
In particular, the following will be used at several points:
\smallskip

\noindent\textbf{Precompactness criterion.}\enspace \textsl{If $H\subset\Iso(M,\mathbf g)$ is a subgroup
that has one orbit in $\mathcal F(M)$ which is contained in a
compact subset of $\mathcal F(M)$, then $H$ is precompact}. For instance, if $H$ preserves some Riemannian
metric on $M$ and it leaves a non empty compact subset of $M$ invariant, then $H$ is precompact.
\begin{lem}\label{thm:elem}
Let $(M,\mathbf g)$ be a compact Lorentzian manifold and $K$ be a Killing vector field on $M$.
If $K$ is timelike at some point, then it generates a precompact $1$-parameter subgroup of isometries in
$\Iso_0(M,\mathbf g)$.
\end{lem}
\begin{proof}
Let $p\in M$ be such that $\mathbf g\big(K(p),K(p)\big)<0$. Consider the compact
subsets of $TM$ given by:
\[V=\Big\{K(q):q\in M\ \text{is such that}\ \mathbf g\big(K(q),K(q)\big)=\mathbf g\big(K(p),K(p)\big)\Big\},\]
and
\[V^\perp\!=\!\Big\{\!v\in K(q)^\perp:q\in M\ \text{s.t.}\  \mathbf g\big(K(q),K(q)\big)\!=\!
\mathbf g\big(K(p),K(p)\big),\ \mathbf g(v,v)=1\Big\}.\]
Consider an orthogonal basis $b=(v_1,\ldots,v_n)$ of $T_pM$ with $v_1=K(p)$ and $\mathbf g(v_i,v_j)=\delta_{ij}$ for
$i,j\in\{2,\ldots,n\}$. The $1$-parameter subgroup generated by $K$ in $\Iso(M,\mathbf g)$ can be identified with the
$\R$-orbit of the basis $b$ by the action of the flow of $K$ on the frame bundle $\mathcal F(M)$. Every vector of a
basis of the orbit belongs to
the compact subset $V\bigcup V^\perp$, and this implies that the orbit of $b$ is precompact in the frame
bundle $\mathcal F(M)$.
\end{proof}
\subsection{An algebraic criterion for precompactness}
Now observe that if a compact manifold $(M,\mathbf g)$ admits a Killing vector field which is timelike somewhere,
then, by continuity, sufficiently close Killing fields are also timelike somewhere.
Thus, if one wants to study the (connected component of the) isometry group of a
Lorentz manifold that has a Killing vector field which is timelike at some point,
it is a natural question to ask is which (connected) Lie groups have open sets of
precompact 1-parameter subgroups. The problem is better cast in terms of the Lie
algebra; we will settle this question in our next:
\begin{prop}\label{thm:charprecompactsubgr}
Let $G$ be a connected Lie group, $K\subset G$ be a maximal compact subgroup, and
$\mathfrak k\subset\mathfrak g$ be their Lie algebras. Let $\mathfrak m$ be a $\mathrm{Ad}_K$-invariant
complement of $\mathfrak k$ in $\mathfrak g$.
Then, $\mathfrak g$ has a non empty open cone of vectors that generate precompact $1$-parameter
subgroups of $G$ if and only if there exists $\mathfrak v\in\mathfrak k$ such that
the restriction of $\mathrm{ad}_{\mathfrak v}:\mathfrak m\to\mathfrak m$ is an isomorphism.
\end{prop}
\begin{proof}
Let $\mathfrak C\subset\mathfrak g$ be the cone of vectors that generate precompact $1$-parameter subgroups
of $G$; we want to know when $\mathfrak C$ has non empty interior. Clearly $\mathfrak C$ contains $\mathfrak k$, and
every
element of $\mathfrak C$ is contained in the Lie algebra $\mathfrak k'$ of some maximal compact subgroup $K'$ of $G$.
Since all maximal compact subgroups of $G$ are conjugate (see for instance \cite{Hoch}), it follows that $\mathfrak
C=\mathrm{Ad}_G(\mathfrak k)$,
i.e., $\mathfrak C$ is the image of the map $F:G\times\mathfrak k\to\mathfrak g$ given by $F(g,\mathfrak
v)=\mathrm{Ad}_g(\mathfrak v)$.
We claim that $\mathfrak C$ has non empty interior if and only if the differential $\mathrm dF$ has maximal rank at
some
point $(g,\mathfrak v)\in G\times\mathfrak k$. The condition is clearly sufficient, and by Sard's theorem is
also necessary; namely, if $\mathrm dF$ has never maximal rank then all the values of $F$ are critical, and
they must form a set with empty interior.
The second claim is that it suffices to look at the rank of $\mathrm dF$ at the points $(e,\mathfrak v)$,
where $e$ is the identity of $G$.
This follows easily observing that the function is $G$-equivariant in the first variable.
Now, the differential of $F$ at $(e,\mathfrak v)$ is easily computed as:
\[\mathrm dF_{(e,\mathfrak v)}(\mathfrak g,\mathfrak k)=[\mathfrak g,\mathfrak v]+\mathfrak k=
[\mathfrak m,\mathfrak v]+\mathfrak k.\]
Thus, $\mathrm dF_{(e,\mathfrak v)}$ is surjective if and only if there exists $\mathfrak v\in\mathfrak k$ such
that $[\mathfrak m,\mathfrak v]=\mathfrak m$, which concludes the proof.
\end{proof}
\end{section}
\begin{section}{The identity connected component of the isometry group}
\label{sec:identitycomponent}
The geometric structure of compact Lorentz manifold whose isometry group contains a group which
is locally isomorphic to an \emph{oscillator group} or to $\SL(2,\R)$ is well known. Let us recall
(see \cite[\S 1.6]{Ghani1}) that a compact Lorentz manifold that admits a faithful isometric action
of a group locally isomorphic to $\SL(2,\R)$ has universal cover which is given by a warped product
of the universal cover of $\SL(2,\R)$, endowed with the bi-invariant Lorentz metric given by
its Killing form, and a Riemannian manifold. Every such manifold admits everywhere timelike Killing vector field,
corresponding to the timelike vectors of the Lie algebra $\mathfrak{sl}(2,\R)$.

Oscillator groups are characterized as the only simply-connected solvable and non commutative Lie groups
that admit a bi-invariant Lorentz metric (see \cite{MedRev}); oscillator groups possess a lattice,
i.e., a co-compact discrete subgroup.
More precisely, an oscillator group $G$ is a semi-direct product
$\mathds S^1\ltimes\mathrm{Heis}$, where $\mathrm{Heis}$ is a Heisenberg group  (of some dimension $2d +1$). There
are
positivity condition on the eigenvalues of the automorphic $\mathds S^1$ action on the Lie algebra
$\mathfrak{heis}$ (ensuring the existence of a bi-invariant Lorentz metric), and arithmetic
conditions on them (ensuring existence of a lattice).

It is interesting and useful to consider oscillator groups as objects completely similar to
$\SL(2,\R)$, from a Lorentz geometry viewpoint.  In particular, regarding our arguments in the present
paper, both cases are perfectly parallel. Let us notice however, some differences (but with no incidence
on our investigation here). First, of course, the bi-invariant Lorentz metrics on an oscillator group do
not correspond to its Killing form, since this latter is degenerate (because the group is solvable).
Another fact, in non-uniqueness of these bi-invariant metrics, but surprisingly, their uniqueness
up to automorphism. In the $\SL(2,\R)$-case, we have uniqueness up to a multiplicative constant.
Also, we have essential uniqueness of lattices in an oscillator group, versus their abundance in $\SL(2,\R)$.

Let us describe briefly the construction of Lorentz manifolds endowed with a faithful isometric $G$-action,
where $G$ is either $\SL(2,\R)$ or an oscillator group. The construction
starts by considering right quotients $G/\Gamma$, where $\Gamma$ is a lattice of
$G$. The $G$-left action is isometric exactly  because the metric is bi-invariant.
A slight generalization is obtained by considering a Riemannian manifold $(\widetilde N,\widetilde{\mathbf g})$
and quotients of the direct metric product $X= \widetilde{N} \times  G$ by a
discrete subgroup $\Gamma$ of $\mathrm{Iso}(\widetilde N,\widetilde{\mathbf g})\times G$.
Observe here that since the isometry group of the Lorentz manifold $X$
is $\mathrm{Iso}(\widetilde N,\widetilde{\mathbf g})\times (G \times G)$,
it is possible to take a quotient by a subgroup $\Gamma$ contained in this full group.
The point is that we assumed $G$ acting (on the left) on the quotient, and hence,
$G$ normalizes $\Gamma$; but since $G$ is connected, it centralizes $\Gamma$. Therefore, only the
right $G$ factor in the full group remains (since the centralizer of the left action is exactly the
right factor). Observe however that $\Gamma$ does not necessarily split.
Indeed, there are examples where $\Gamma$ is discrete co-compact in
$\mathrm{Iso}(\widetilde N,\widetilde{\mathbf g})\times G$, but its projection on each factor is dense!

Next, warped products yield a more general construction. Rather than a direct product metric
$\widetilde{\mathbf g}\oplus\kappa$, one endows $\widetilde{N}\times G$ with
a metric of the form $\widetilde{\mathbf g}\oplus w \kappa$, where $w$ is a positive function on
$\widetilde{N}$, and  $\kappa$ is the bi-invariant Lorentz metric on $G$.
Here, there is one difference between the case of $\SL(2, \R)$ and the oscillator case.
For $\SL(2,\R)$ this is the more general construction, but in the oscillator case, some ``mixing''
between $G$ and $\widetilde{N}$ and  also a mixing of their metrics, is also possible, see \cite[\S 1.2]{Ghani1}.
Finally,  to be more accurate, it should be emphasized that it is the
universal cover of $G$ that must be considered in these constructions.
(It turns out however that, it is a finite cover of $G$ which acts faithfully and
not the universal cover. This is because of the compactness assumption on the Lorentz manifold $M$, see \cite[\S
2]{Ghani1}).
\smallskip

Let us study now the situation when the isometry group does not contain any group which is locally
isomorphic to $\SL(2,\R)$ or to an oscillatory group.
\begin{theo}\label{theorem1} Let $(M,\mathbf g)$ be a  compact Lorentz manifold that admits a Killing vector field
which is
timelike at some point.  Then, the identity component of its isometry group is compact,
unless it contains a group locally isomorphic to $\SL(2, \R)$ or to an oscillator group.
\end{theo}
\begin{proof}
By the classification result in \cite{AdamsStuck, Ghani1}, if $\Iso_0(M,\mathbf g)$ does not contain a group locally
isomorphic to $\SL(2,\R)$ or to an oscillator group, then $\mathfrak{Iso}(M,\mathbf g)$ can be written
as a Lie algebra direct sum $\mathfrak h+\mathfrak a+\mathfrak c$, where $\mathfrak h$ is a Heisenberg
algebra, $\mathfrak a$ is abelian, and $\mathfrak c$ is semi-simple and compact. Our aim is to show that
the Heisenberg summand $\mathfrak h$ in fact does not occur in the decomposition, and that the abelian
group $A$ corresponding to the summand $\mathfrak a$ is compact. By the assumption that
$(M,\mathbf g)$ has a Killing vector field which is timelike at some point, $\mathfrak{Iso}(M,\mathbf g)$
must contain a non empty open cone of vectors that generate a precompact $1$-parameter subgroup
of $\Iso_0(M,\mathbf g)$ (Lemma~\ref{thm:elem}). The first observation is that, since $\mathfrak c$ is compact, if
$\mathfrak h+\mathfrak a+\mathfrak c$ has an
open cone of vectors that generate precompact $1$-parameter subgroups, than so does the subalgebra
$\mathfrak h+\mathfrak a$. Moreover, by the same compactness argument, we can also assume that
the abelian Lie subgroup $A$ is simply connected. The proof of our result will be concluded once we
show that given any Lie group $G$ with Lie algebra $\mathfrak g=\mathfrak h+\mathfrak a$, $\mathfrak h$ and
$\mathfrak a$ as above, does not have an open set of precompact $1$-parameter subgroups. To this aim,
write $\mathfrak h+\mathfrak a=\mathfrak m+\mathfrak k$, with $\mathfrak k$ the Lie algebra of a maximal compact
subgroup
$K$ of $G$ and $\mathfrak m$ a $\mathfrak k$-invariant complement of $\mathfrak k$ in $\mathfrak g$.
If either $\mathfrak h$ or $\mathfrak a$ is not zero, then also $\mathfrak m$ is non zero.
Since $\mathfrak h+\mathfrak a$ is nilpotent,  for no $\mathfrak x\in\mathfrak k$ the map $\mathrm{ad}_{\mathfrak
x}:\mathfrak m\to\mathfrak m$ is injective, and by Proposition~\ref{thm:charprecompactsubgr}, $G$ does
not have an open set of precompact $1$-parameter subgroups.
\end{proof}

\end{section}
\begin{section}{When the isometry group has infinitely many connected components}
Let us now study the situation when the isometry group of a Lorentzian manifold with a timelike
Klling vector field has infinitely many connected components.
\begin{exa}\label{rem:exampleflattorus}
An important example is that of a flat torus, see for instance \cite{ThierryGhani}.
Consider a Lorentzian scalar product $\mathbf g_0$ in $\R^n$, and let the quotient
$M=\mathds T^n=\mathds R^n/\mathds Z^n$ be endowed with the induced metric, still denoted by $\mathbf g_0$.
Then, $\mathrm{Iso}(M,\mathbf g_0)$ is isomorphic to the semi-direct product $\mathds T^n\rtimes\mathrm O(\mathbf
g_0,\mathds Z)$,
where $\mathrm O(\mathbf g_0,\mathds Z)$ is the infinite discrete group $\mathrm O(\mathbf
g_0)\cap\mathrm{GL}(n,\mathds Z)$.
\end{exa}
\subsection{Compactness of the identity connected component}
A first consequence of Theorem~\ref{theorem1} is that when $\Iso(M,\mathbf g)$ has infinitely many
connected components, then the connected component of the identity has to be compact. Before proving
this, let us recall a general fact on the isometry group of a covering.
If $(M,\mathbf g)$ is a pseudo-Riemannian manifold and $\pi:\widetilde M\to M$ is a covering, let
$\widetilde{\mathbf g}$ be the pseudo-Riemannian metric on $\widetilde M$ given by the pull-back of $\mathbf g$ by
$\pi$, so that $\pi$ becomes a local isometry. Let $H\subset\Iso(\widetilde M,\widetilde{\mathbf g})$ be the closed
subgroup consisting of all isometries
$\widetilde f:\widetilde M\to\widetilde M$ that \emph{descend to the quotient $M$}, i.e.,
such that $\pi\big(\tilde f(x)\big)=\pi\big(\tilde f(y)\big)$ whenever $\pi(x)=\pi(y)$.
There is a natural homomorphism $\varphi:H\to\Iso(M,\mathbf g)$ that associates to each $\widetilde f\in H$ the
unique isometry $f$ of $M$ such that $f\circ\pi=\pi\circ\widetilde f$. The kernel of $\varphi$ is
a discrete subgroup $\Gamma$ of $H$ given by the group of covering automorphisms of $\pi$.
If $\pi$ is the universal covering of $M$, then the homomorphism $\varphi$ is surjective, and
$H$ is precisely the \emph{normalizer} $\Nor(\Gamma)$ of $\Gamma$ in $\Iso(\widetilde M,\widetilde{\mathbf g})$,
so that $\Iso(M,\mathbf g)$ is isomorphic to the quotient $\Nor(\Gamma)/\Gamma$.
\begin{prop}\label{thm:propinfmanyconncomp} Let $(M,\mathbf g)$ be a compact Lorentzian manifold that admits a
somewhere timelike Killing vector field.
If $\Iso(M,\mathbf g)$ has infinitely many connected components, then $\Iso_0(M,\mathbf g)$ is compact.
\end{prop}
\begin{proof}
By Theorem~\ref{theorem1}, if $\Iso_0(M,\mathbf g)$ is not compact, then it contains a group which is locally
isomorphic
to $\SL(2,\R)$ or to an oscillator group, in which case the geometric structure of $(M,\mathbf g)$ is well understood
(see \cite[Theorem~1.13, Theorem~1.14]{Ghani1}).


Assume that $\Iso_0(M,\mathbf g)$ contains a group locally isomorphic to $G= \SL(2,\R)$.
Then, there exists a Riemannian manifold $(\widetilde N,\widetilde{\mathbf g})$
such that, up to a finite cover, $M$ is a quotient of $X = \widetilde{N} \times G$ by a subgroup  $\Gamma$ of
 $\Iso(\widetilde N,\widetilde{\mathbf g})\times G$. Here, $G$ is endowed with the bi-invariant Lorentzian metric
 defined by the Killing form of $\mathrm{SL}(2,\R)$ and the metric of $\widetilde{N}\times G$ may be warped
 rather than a direct product. Any isometry of $X$ that normalizes the $G$-action on $X$
 is the product of an isometry $a$ of the first factor $\widetilde N$, and an isometry
 of $G$ of the form $(b,c)\in G\times G$, with $b$ acting by right multiplication and $c$ by left multiplication
 on $G$. When the metric of $X$ is warped, then $a$ must preserve the warping function.
 We will denote by $(a,b,c)$ such an isometry of $X$.

 The fundamental group $\Gamma$ commutes with the $G$-action, and therefore its elements have the form
 $(a,b,1)$, for they must commute with all elements of the form $(1,1,c)$.
 One can define a \emph{Riemannian} metric on $X$ by replacing the Lorentz metric of the factor $G$
 with a right invariant Riemannian metric on $G$, that will be denoted by $\mathbf m$.
 By construction, such a new Riemannian  metric on $X$ is preserved by $\Gamma$, thus it descends
 to a Riemannian metric $\mathbf h$ on $M$. Now, the map $(a,b,c)\mapsto(a,b)$ defines a homomorphism
 \[\phi:\mathrm{Iso}(M,\mathbf g)\longrightarrow\mathrm{Iso}(M,\mathbf h),\]
 whose kernel is exactly $G$. Observe however that $\mathrm{Iso}(M,\mathbf h)$ is compact, but
 $\phi$ may not be onto. Denote by $H_1$ the closed subgroup of $\mathrm{Iso}(M,\mathbf h)$
 consisting of elements that preserve the foliation of $M$ determined by $G$. This is a closed
 subgroup, and all its elements have the split form $(a,d)$, where $d\in\mathrm{Iso}(G,\mathbf m)$.

 Now, $\mathrm{Iso}(G,\mathbf m)$ contains $G$ (acting by right multiplication), but it may be strictly
 larger than $G$ when the isotropy group of $(G,\mathbf m)$ is non trivial, and a further reduction is needed.
 More precisely, let $e_1,e_2,e_3$ be a basis of the Lie algebra of $G$, and let $E_1,E_2,E_3$ be the corresponding
 vector fields on $M$. Denote by $H_2$ the closed subgroup of $H_1$ consisting of isometries $f\in\mathrm{Iso}(M,\mathbf h)$ such that $f_*E_i=E_i$, $i=1,2,3$. The lift to $X$ of such an isometry $f$
 has the form $(a,b,1)$, thus we have an exact sequence $1\to G\to\mathrm{Iso}(M,\mathbf g)\to H_2\to 1$.
 Since $H_2$ is compact, it has a finite number of connected components. Moreover, $G$ is connected,
 it follows that $\mathrm{Iso}(M,\mathbf g)$ has the same number of connected components as $H_2$, and
 we are done.

The case where $\Iso_0(M,\mathbf g)$ contains an oscillator group is treated following
the same arguments.
\end{proof}
\newpage

\subsection{A further compactness result}
\begin{lem}\label{thm:Riemmetrh0}
Let $G$ be a Lie group acting on a manifold $R$, and let $G_0$ be a compact normal subgroup
of $G$ all of whose orbits in $R$ have the same dimension. Then:
\begin{itemize}
\item the distribution $\Delta$ tangent to the $G_0$-orbits is smooth, and it is preserved by
$G$;
\item there exists a Riemannian metric $\mathbf h_0$ on $\Delta$  which is preserved by $G_0$ and by the
    centralizer $\Centr(G_0)$ of $G_0$ in $G$.
\end{itemize}
\end{lem}
\begin{proof}
Introduce the following notation: for $x\in R$, let $\beta_x:G\to R$ be the map $\beta_x(g)=g\cdot x$, and
let $L_x:\mathfrak g\to T_xM$ be its differential at the identity. Here $\mathfrak g$ is the Lie algebra
of $G$. The map $M\times\mathfrak g\ni(x,\mathfrak v)\mapsto L_x(\mathfrak v)\in TM$ is a smooth vector
bundle morphism from the trivial bundle $M\times\mathfrak g$ to $TM$. The distribution $\Delta$ is the image
of the sub-bundle $M\times\mathfrak g_0$, where $\mathfrak g_0\subset\mathfrak g$ is the Lie algebra of $G_0$.
Since the orbits of $G_0$ have the same dimension, then the image of $M\times\mathfrak g_0$
is a smooth sub-bundle of $TM$ (recall that the image of a vector bundle morphism is a smooth sub-bundle if it has
constant rank). The action of $G$ preserves $\Delta$ because $G_0$ is normal, which concludes the proof
 of the first assertion.

The construction of $\mathbf h_0$ goes as follows. Choose a positive definite inner product $B$
on $\mathfrak g_0$ which is $\mathrm{Ad}_{G_0}$-invariant; the existence of such $B$ follows from
the compactness of $G_0$. For all $x\in R$, the restriction to $\mathfrak g_0$ of $L_x$ gives a
surjection $L_x\vert_{\mathfrak g_0}:\mathfrak g_0\to \Delta_x$; denote by $V_x$ the $B$-orthogonal
complement of the kernel of this map, given by $\Ker(L_x\vert_{\mathfrak g_0})=\Ker(L_x)\cap\mathfrak g_0$.
The value of $\mathbf h_0$ on $\Delta_x$ is defined
to be the push-forward via the map $L_x$ of the restriction of $B$ to $V_x$.
In order to see that such metric is invariant by the action of $G_0$ and of its centralizer,
for $g\in G$ denote by $\mathcal I_g:G_0\to G_0$ the conjugation by $g$ (recall that $G_0$ is normal)
and by $\gamma_g:M\to M$ the diffeomorphism $x\mapsto g\cdot x$; for fixed $x\in M$ we have a commutative
diagram:
\[\xymatrix{G_0\ar[r]^{\mathcal I_g}\ar[d]_{\beta_x}&G_0\ar[d]^{\beta_{gx}}\cr M\ar[r]_{\gamma_g}&M}\]
Differentiating at the identity the diagram above we get:
\[\xymatrix{\mathfrak g_0\ar[rr]^{\mathrm{Ad}_g}\ar[d]_{L_x}&&\mathfrak g_0\ar[d]^{L_{g\cdot x}}\cr
\Delta_x\ar[rr]_{\mathrm d\gamma_g(x)}&&\Delta_{g\cdot x}}\]
Now assume that $g$ is such that $\mathrm{Ad}_g$ preserves $B$; this holds by assumption when
$g\in G_0$, and clearly also for $g$ in the centralizer of $G_0$ (in which case $\mathrm{Ad}_g$ is the
identity!). For such a $g$, since $\mathrm{Ad}_g\big(\Ker(L_x)\cap\mathfrak g_0\big)=\Ker(L_{g\cdot x})$,
then also $\mathrm{Ad}_g(V_x)=V_{g\cdot x}$, and thus we have a commutative diagram:
\[\xymatrix{V_x\ar[rr]^{\mathrm{Ad}_g}\ar[d]_{L_x}&& V_{g\cdot x}\ar[d]^{L_{g\cdot x}}\cr
\Delta_x\ar[rr]_{\mathrm d\gamma_g(x)}&&\Delta_{g\cdot x}},\]
from which it follows that $\mathrm d\gamma_g(x)$ preserves the metric $\mathbf h_0$, proving
the last statement in the thesis.
\end{proof}

\noindent
Let us now assume throughout that $(M,\mathbf g)$ is a compact Lorentzian manifold
that admits a Killing vector field $K$ which is timelike somewhere, and with $\Iso_0(M,\mathbf g)$ compact.
To simplify notations, let us set $G=\Iso(M,\mathbf g)$, $G_0=\Iso_0(M,\mathbf g)$, $\mathfrak
g=\mathfrak{Iso}(M,g)$,
$\mathrm{Aut}(G_0)$ the group of automorphisms
of $G_0$, $\mathrm{Inn}(G_0)$ the normal subgroup of inner automorphisms of $G_0$, and
$\mathrm{Out}(G_0)=\mathrm{Aut}(G_0)/\mathrm{Inn}(G_0)$. Since $G_0$ is normal in $G$, then $G$ acts
 by conjugation on $G_0$, and we have a homomorphism $G\to\mathrm{Aut}(G_0)$, and composing with the
 quotient map we have a homomorphism:
\begin{equation}\label{eq:defrho}
\rho:G\longrightarrow\mathrm{Out}(G_0),
\end{equation}
whose kernel $G_1=\mathrm{Ker}(\rho)$ is the product $G_0\cdot\Centr(G_0)$.

\begin{prop}\label{thm:G1compact}
$G_1$ is compact.
\end{prop}
\begin{proof}
Let $R$ be the (non empty) open subset of $M$ consisting of all points $x$ whose $G_0$-orbit
$\mathcal O(x)$ has maximal dimension (among all $G_0$-orbits), and such that $\mathcal O(x)$
is timelike, i.e., the restriction of $\mathbf g$ to $\mathcal O(x)$ is Lorentzian.
Recall that the set of points whose $G_0$-orbit has maximal dimension
is open and dense, and $R$ is the intersection
of this dense open subset with the open subset of $M$ where $K$ is timelike.
We claim that there exists a Riemannian metric $\mathbf h$ on $R$ which is preserved by $G_1$.
Such a metric $\mathbf h$ is constructed as follows: on the distribution $\Delta$ tangent to the
$G_0$-orbits is given by the metric $\mathbf h_0$ as in Lemma~\ref{thm:Riemmetrh0}, on the
$\mathbf g$-orthogonal distribution $\Delta'$ it is the (positive definite) restriction of $\mathbf g$,
furthermore $\Delta$ and $\Delta'$ are declared to be $\mathbf h$-orthogonal.
Note that the distributions $\Delta$ and $\Delta'$ are preserved by $G$, the metric $\mathbf g$ on $\Delta'$
is preserved by all elements of $g$, and the metric $\mathbf h_0$ on $\Delta$ is preserved by
all elements of $G_1$, which proves our claim.

Observe that $G_1$ is closed in $G$.
Using the precompactness criterion, in order to prove the compactness of $G_1$ it
suffices to show that $G_1$ leaves some compact subset of $R$ invariant. If $d$ is the dimension
of the principal $G_0$-orbits, let $\mathrm{vol}_d$ denote $d$-dimensional volume induced by
the Lorentzian restriction of $\mathbf g$ on the $G_0$-orbits in $R$. Set $a=\sup\limits_{x\in
M}\mathrm{vol}_d\big(\mathcal O(x)\big)$, and let $T$ the compact subset of $M$ consisting of all points $x$ such
that
$\mathrm{vol}_d\big(\mathcal O(x)\big)=a$. Note that the volume function is lower semi-continuous,
see Lemma~\ref{thm:interestinitsown} below,
and it admits maximum in the compact set of all $G_0$-orbits, so that $a$ is well defined and $T$ is not empty.
Moreover, since the condition $\mathrm{vol}_d\big(\mathcal O(x)\big)\ge a$ is closed, then $T$ is compact.
Clearly, $T$ is contained in $R$, and it is preserved by $G_1$. This concludes the proof.
\end{proof}
We have used the following result, which has some interest in its own:
\begin{lem}\label{thm:interestinitsown}
Let $G$ be a compact Lie group that acts smoothly by isometries on a Lorentzian manifold $(L,\mathbf g_L)$ in such
a way that:
\begin{itemize}
\item[(a)] all the $G$-orbits have the same dimension;
\item[(b)] every $G$-orbit is timelike, i.e., the restriction of the metric $\mathbf g_L$ to each orbit
is Lorentzian.
\end{itemize}
Then, the function $L\ni x\mapsto\mathrm{vol}_x\in\R$ is lower semi-continuous, where
$\mathrm{vol}_x$ denotes the volume of the orbit $\mathcal O(x)$ relative to the volume induced by the restriction
to $\mathcal O(x)$ of the Lorentzian metric $\mathbf g_L$.
\end{lem}
\begin{proof}
By standard result on group actions (see for instance \cite{Bredon}), through every $x\in L$ there is a \emph{slice}
for the action of $G$, i.e., a smooth submanifold $S_x\subset L$ containing $x$ that has (among others) the following
properties:
\begin{enumerate}
\item\label{itm:lemn1} the map $G\times S_x\ni(g,y)\mapsto g\cdot y\in L$ is open;
\item\label{itm:lemn2} denoting by $G_x$ the isotropy of $x$, if $g\not\in G_x$ then $g\cdot S_x\cap
    S_x=\emptyset$.
\end{enumerate}
By \eqref{itm:lemn1}, if $V\subset S_x$ is a neighborhood of $x$ in $S_x$, then $G\cdot V$ is a neighborhood
of $x$ in $L$; thus, since the function volume is constant on each orbit, it suffices to show that the map $S_x\ni
y\mapsto\mathrm{vol}_y\in\R$ is lower
semi-continuous.

By \eqref{itm:lemn2}, if $y\in S_x$, then $G_y\subset G_x$. Let $H\subset G_x$ be the identity connected component
of $G_x$. By assumption (a), $\mathrm{dim}\big(\mathcal O(y)\big)=\mathrm{dim}\big(\mathcal O(x)\big)$, from
which it follows that $H$ is contained in $G_y$. We therefore have a well defined map
$\overline\beta:G/H\times S_x\to L$, given by $\overline\beta(gH,y)=g\cdot y$. This map is smooth, because
$\mathfrak q:G\times S_x\to (G/H)\times S_x$ is a smooth surjective submersion, and
$\overline\beta\circ\mathfrak q:G\times S_x\to L$ is smooth.

Fr $y\in S_x$, consider the smooth map $\overline\beta_y:G/H\to L$ defined by
$\overline\beta_y(gH)=\overline\beta(gH,y)$;
its image is given by the embedded submanifold $\mathcal O(y)$, thus we have a smooth map
$\overline\beta_y:G/H\to\mathcal O(y)$. Such map is a finite covering, with folding (cardinality of the fiber)
equal to $\big\vert G_y/H\big\vert$. Namely, $\overline\beta_y$ can be written as the composition of the
diffeomorphism
$G/G_y\ni gG_y\mapsto g\cdot y\in\mathcal O(y)$ with the covering map $G/H\to G/G_y$.
The latter is a smooth fibration, with discrete fiber, hence a covering map, with folding equal
to $\big\vert G_y/H\big\vert$ (see for instance \cite[Proposition~2.1.14]{MaslovBook}).

It follows that the map $\overline\beta_y:G/H\to L$ is an immersion (composition of an immersion and a local
diffeomorphism),
and we obtain a symmetric $(0,2)$-tensor on $G/H$ given by the pull-back $\mathbf g_y=\overline\beta_y^*(\mathbf
g_L)$.
In fact, we have a smoothly varying family $S_x\ni y\mapsto\mathbf g_y$ of symmetric $(0,2)$-tensors on $G/H$.
Such tensors are nondegenerate and Lorentzian, in fact $\mathbf g_y$ is also given as the pull-back
by the covering map $\overline\beta_y:G/H\to\mathcal O(y)$ of the restriction to $\mathcal O(y)$ of $\mathbf g_L$.
It follows that the volume of $G/H$ relatively to the metric $\mathbf g_y$, which is a continuous function
of $y$ (volume relative to a continuous family of measures), is equal to $\big\vert G_y/H\big\vert\mathrm{vol}_y$.

The desired semi-continuity property of the function $y\mapsto\mathrm{vol}_y$ follows now easily observing
that for $y\in S_x$, $G_y\subset G_x$, thus $\big\vert G_y/H\big\vert\mathrm{vol}_y\le
\big\vert G_x/H\big\vert\mathrm{vol}_y$.
\end{proof}
\section{The toral factor}
\label{sec:toralfactor}
We will now pursue the study of the action of subgroups of the isometry group of $M$.
The idea is to consider a suitable \emph{reduction} of $\mathrm{Iso}(M,\mathbf g)$, i.e., a closed subgroup
$G\subset\mathrm{Iso}(M,\mathbf g)$ such that, denoting by $G_0$ its identity connected component,
the following properties are satisfied:
\begin{eqnarray}
\notag
&&\text{$G$ is non compact,}\\ \label{eq:situation}
&&\text{$G_0$ is compact,}\\ \notag
&&\text{$G_0$ has some timelike orbit.}
\end{eqnarray}
We observe that \eqref{eq:situation} is preserved by passing to finite index subgroups of $G$.
Namely, if $G'\subset G$ has finite index, then it is closed, non compact, and it has the same
connected component of the identity as $G$.
\subsection{Reduction of $G^0$}
Recall the homomorphism $\rho:G\to\mathrm{Out}(G_0)$ in \eqref{eq:defrho}; by  Proposition~\ref{thm:G1compact},
$G_1=\mathrm{Ker}(\rho)$ is compact, and in particular $\rho(G)$ is non compact.

If $G_0$ has an almost decomposition $\mathds T^k\times K$, where $K$ is semisimple, then the image of $\rho$ is
contained in $\mathrm{Out}(\T^k)=\GL(k,\Z)$.

Since $K$ is normal, we have a representation $r:G\to\mathrm{Aut}(K)$; let  $G^\prime$ be its kernel, which
is the centralizer of $K$ in $G$. One can see that $G/G^\prime= K/\mathrm{Centr}(K)$. Indeed, since
$\mathrm{Aut}(K)=\mathrm{Int}(K)$, for any $f\in G$, there exists $k\in K$ such that $r(f)=r(k)$,
that is, $fk^{-1}\in G^\prime$. Therefore, $G/ G^\prime$ is a quotient of $K$.
This quotient is easily identified with $K /\mathrm{Centr}(K)$.

In some sense, going from $G$ to $G^\prime$ allows one to kill the semi-simple factor, that is to assume that
the identity component is a torus, and that the discrete part $G/G^0$ has not changed.
More precisely,
let us now describe how to ``forget'' the semisimple factor $K$ keeping the identity component with somewhere
timelike orbits. Let $X$ be a somewhere timelike Killing field.
The closure of its flow is a product (possibly trivial) of  two tori, $K_1 \times K_2$, where $K_1 $ (resp., $K_2$)
is a subgroup of $\T^k$ (resp., of $K$). Since $G^\prime$ centralizes $K_2$, we have a  direct product
group $G^\prime \times K_2$.

Summarizing, we have proven the following:

\begin{lem} There  is a subgroup $G$ of  $\mathrm{Iso}(M,\mathbf g)$ satisfying \eqref{eq:situation}
and having an abelian identity component $G_0=\T^k$.
\end{lem}

With such reduction of the group $G$, we can now consider the action of $G$ on $G_0\cong\mathds T^k$ given by the
representation $\rho:G\to\mathrm{Out}(\mathds T^k)=\mathrm{GL}(k,\Z)$; in order to distinguish
the action of $G$ on $M$ and on $\mathds T^k$, we will call the latter the \emph{$\rho$-action}.
\begin{cor}\label{thm:notorsion}
Up to a finite index reduction, the quotient group $\Gamma=G/G_0$ is torsion free, i.e., all its non trivial
elements have infinite order.
\end{cor}
\begin{proof}
Choose any torsion free finite index subgroup $H$ of $\mathrm{GL}(k,\Z)$ (it exists by Selberg lemma \cite{Roger}),
and set $G'=\rho^{-1}(H)$. This is a finite index torsion-free subgroup of $G$.
\end{proof}

\subsection{Preliminary properties of the  $\T^k$-action}
\begin{lem} The $\T^k$-action on $M$ is locally free on a dense open set.
\end{lem}
\begin{proof}  It follows immediately from Corollary~\ref{thm:faithfulfreeaction}.
\end{proof}

\begin{lem}\label{thm:preserveLorntzform}
After replacing $G$ by a finite index subgroup, the $\rho$-action on $\T^k$
preserves some Lorentz metric. In particular, one can see $\rho(G)$ as lying in $\GL(k,\Z)\cap\mathrm{SO}(q)$,
where $q$ is a Lorentz form on $\R^k$.
\end{lem}
\begin{proof} The $\rho$-action of $G$ on $G_0\cong\T^k$ by conjugation induces an action of $G$ on the
space $\mathrm{sym}(\R^k)$ of symmetric bilinear forms on $\mathds R^k$. By \eqref{eq:equivGauss},
the compact subset given by the image of the Gauss map $\mathcal G$ is invariant by this action;
moreover, the image of the Lebesgue measure of $M$ by the Gauss map is also invariant by $\rho(G)$.
Such compact subset contains a non empty open subset consisting of Lorentz forms,
because $G_0$ has timelike orbits in $M$.
By Corollary~\ref{thm:finiteindexsubgroup}, there exists a finite index subgroup of $H'$ of $\rho(G)$
that preserves some Lorentz form on $\R^k$. The desired finite index subgroup of $G$ is $\rho^{-1}(H')$.
\end{proof}

\section{Actions of almost cyclic groups}
Choose $f\in G$, $f\not\in G_0$; then, $\rho(f)\in\mathrm{GL}(k,\Z)$ has infinite order by
Corollary~\ref{thm:notorsion}.
Consider the group $G = G_f$ generated by $f$ and $\T^k$.
Up to a compact normal subgroup, $G_f$ is cyclic, which justifies the name \emph{almost cyclic}.
One can  prove:
\begin{lem}  $G=G_f$ is a closed subgroup of $\mathrm{Iso}(M,\mathbf g)$ satisfying \eqref{eq:situation}.
It is isomorphic to a semi-direct product $\Z\ltimes \T^k$.
\end{lem}

\begin{lem} \label{hyperbolic.parabolic} If $A\in\mathrm{GL}(k,\Z)$  is parabolic,
then some power of $A$ is rationally equivalent to:
\[\begin{pmatrix}
\begin{matrix}1 &  t & -t^2/2\\
 0 & 1 & t\\
 0 & 0&  1\end{matrix}&\mathbf0 \\
\mathbf0 & \mathrm{Id}_{k-3}  \\
 \end{pmatrix}\]
This means that the subspaces $\{e_1\}$, $\{e_1, e_2\}$, $\{e_1, e_2, e_3\}$ and
$\{e_4, \ldots, e_k\}$ are rational.

In particular, there is an $A$-invariant rational $3$-space, on whose orthogonal, which is not necessarily rational,
the $A$-action is trivial.
\end{lem}
\begin{proof}   The proof is quite standard.  Let $A$ have the normal form above.
Consider $\mathcal E$ the kernel of $(A-1)^3$, which is a rational subspace, and it
contains the subspace $\mathcal E_0 = \{e_1, e_2, e_3\}$.
On $\mathcal E/\mathcal E_0$, $A$ is elliptic, but since it satisfies $(A-1)^3 = 0$, its eigenvalues are roots of 1.
More precisely, replacing $A $ by a $A^3$, we can assume that $A$ is trivial on its elliptic subspace in $\mathcal
E$.

Since $\mathcal E$ is rational, $A$ determines an integer matrix in $\mathrm{GL}(\R^k/\mathcal E)$.
This is an elliptic matrix. So, all the eigenvalues of $A$ are roots of unity, and therefore,
after passing to a power, we can assume that $1$ is the unique eigenvalue of $A$.

Consider $A-1$ and $(A-1)^2$. Their  images are, respectively,
the $2$-plane generated by  $\{e_1, e_3\}$  and the line $\R e_1$.
These two subspaces are thus rational.

The $1$-eigenspace of $A$ is generated by $\{e_1, e_4, \ldots, e_k\}$. It is rational.
We can choose $e_4, \ldots e_k$ rational.
For $e_3$, one can take any rational vector which does not belong to the space
generated by $\{e_1, e_2, e_4, \ldots, e_k\}$.
\end{proof}

\subsection{Structure Theorem}

\begin{theo} \label{almost.cyclic} Let $f\in\mathrm{Iso}(M,\mathbf g)$ act non-periodically on
$\mathrm{Iso}_0(M,\mathbf g)$.
Then, there is a minimal timelike $\rho(f)$-invariant torus $\T^d\subset\mathrm{Iso}_0(M,\mathbf g)$
of dimension $d=3$ or $d \geq 2$ according to whether $\rho(f)$ is  parabolic or hyperbolic, respectively.
The action of $\T^d$ on $M$ is (everywhere) free and timelike.
\end{theo}
We will present the proof of  the theorem in the parabolic case; the hyperbolic case is analogous, in fact,
easier.
So, let $f$ be such that $\rho(f)$ is parabolic. The $3$-torus $\T^d= \T^3$ in
question is the one corresponding to the rational $3$-space associated to $A$ in Lemma~\ref{hyperbolic.parabolic}.
The normal form of $\rho(f)$ on this rational $3$-space is:
\begin{equation}\label{eq:normalformrhof3space}
\rho(f)\cong\begin{pmatrix}1&t&-\tfrac{t^2}2\\0&1&t\\0&0&1\end{pmatrix},\quad t\ne0.
\end{equation}
We need to show that this torus acts freely with timelike orbits on $M$, and the idea is to relate
the dynamics of $f$ on $M$ and the dynamics of $\rho(f)$ on the toral factor.
Towards this goal, we will use the approximately stable foliation of a Lorentz isometry,
introduced in reference~\cite{Zeghib.GAFA1}.
\subsection{Recalls on approximate stability}
Let $\phi$ be a diffeomorphism  of a compact manifold $M$. A vector $v \in T_xM$ is called
\emph{approximately stable} if there is a sequence $v_n \in T_xM$, $v_n \to v$ such that the sequence
$D_x\phi^n v_n$ is bounded in $TM$. The vector $v$ is called \emph{strongly approximately
stable} if $D_x\phi^nv_n \to 0$. The set of approximately stable vectors in
$T_xM$ is denoted $\mathrm{AS}(x, \phi)$, or  sometimes
$\mathrm{AS}(x, \phi, M)$. Their union over $M$ is denoted $\mathrm{AS}(\phi)$, or $\mathrm{AS}(\phi, M)$.
Similarly, $\mathrm{SAS}(x,\phi)$ will denote that set of strongly approximately stable
vectors in $T_xM$, and $\mathrm{SAS}(\phi)=\bigcup_{x\in M}\mathrm{SAS}(x,\phi)$.

The structure of $\mathrm{AS}(\phi)$ when $\phi$ is a Lorentzian isometry has been studied in \cite{Zeghib.GAFA1}:
\begin{theo}[Zeghib \cite{Zeghib.GAFA1}] Let $\phi$ be an isometry of a compact Lorentz manifold $(M,\mathbf g)$
such that the powers $\{\phi^n\}_{n\in\mathds N}$ of $\phi$ form an unbounded set
(i.e., non precompact in $\mathrm{Iso}(M,\mathbf g)$).
Then:
\begin{itemize}
\item $\mathrm{AS}(\phi)$ is a Lipschitz condimension $1$ vector subbundle of $TM$ which is tangent
to a condimension $1$ foliation of $M$ by geodesic lightlike hypersurfaces;
\smallskip

\item $\mathrm{SAS}(\phi)$ is a Lipschitz $1$-dimensional subbundle of $TM$ contained in $\mathrm{AS}(\phi)$
and everywhere lightlike.
\end{itemize}
\end{theo}

\subsection{The action on $M$ vs the toral action}
Denote by $\cT$ the Lie algebra of  $\T^3$, and by $\rho_0(f)$ the linear representation in
associated to $\rho(f)$. More explicitly, $\rho_0(f)$ is the push-forward by $f$ of Killing vector fields,
see formula \eqref{eq:equivariance}.
\begin{lem}   Let $X \in \cT$,  be a Killing field which is approximately stable for
$\rho(f)$ at $1\in\mathds T^3$.   Then,
for all $x \in M$, $X(x) \in T_xM$ is approximately stable. In other words, if
$X\in\mathrm{AS}(0, \rho_0 (f), \cT)$, then $X(x) \in\mathrm{AS}(x, f, M)$ for any $x\in M$.

A totally analogous statement holds for the strong approximate stability.
\end{lem}

\begin{proof} Let  $X_n$  be a sequence of Killing fields in $\cT$ such that $X_n \to X$ and with
$Y_n= f_*^nX_n$ bounded. Clearly $X_n(x)\to X(x)$ for all $x\in M$; moreover, by assumption,
the $Y_n$ are bounded vector fields, and so  $D_x f^n X_n(x) = Y_n(f^nx)$
is bounded, that is $X(x) \in\mathrm{AS}(f)$.
\end{proof}

\begin{lem}\label{thm:Zparabolic}
Assume  $\rho(f)$ parabolic. Then, there is $Z \in \cT$ a Killing field such that
\begin{itemize}
\item[(a)] $Z$ defines a periodic flow $\phi^t$;
\item[(b)] $f$ preserves $Z$, i.e., $f$ commutes with the $1$-parameter group of isometries $\phi^t$
generated by $Z$;
\item[(c)] $Z$ generates the strong approximate stable $1$-dimensional bundle of $f$;
\item[(d)] $Z$ is everywhere isotropic;
\item[(e)] $Z$ is non-singular, hence $Z$ is everywhere lightlike.
\end{itemize}
\end{lem}
\begin{proof} Let $Z$ to be a $1$-eigenvector of $\rho(f)$; since $\rho_0(f) Z=f_*Z= Z$, then
$f$ preserves $Z$.
In the normal form \eqref{eq:normalformrhof3space} of $\rho(f)$,
the vector $Z$ corresponds to the first element of the basis.

The $Z$-direction is rational, since it is the unique $1$-eigendirection of $\rho_0(f)$.
Thus $Z$ defines a periodic flow.

One verifies that $Z$ is strongly approximately stable for the $\rho_0(f)$-action  at $0 \in \cT$. Therefore, at any
$x$ where it does not vanish, $Z(x)$ determines the strongly stable 1-dimensional bundle of $f$. In particular,
$Z(x)$
is isotropic for all $x\in M$. But, non trivial isotropic Killing field cannot
have singularities.
\end{proof}

\subsection{Proof of Theorem~\ref{almost.cyclic} (parabolic case)}
\begin{lem} $\T^3$ preserves the approximate stable foliation $\cF$ of $f$.
\end{lem}

\begin{proof} The group $G_f$ generated by $\mathds T^3$ and $f$ is amenable (it is an extension of the
abelian $\T^3$ by the abelian $\Z$). The statement follows then from \cite[Theorem~2.4, Theorem~2.6]{Zeghib.GAFA1}.
\end{proof}

\begin{lem}
The $\T^3$-action is locally  free.
\end{lem}
\begin{proof}
Let $\Sigma$ be the set of points $x$  having a stabilizer  $S_x$ of  positive dimension.
We claim that if $\Sigma$ is non empty, then there must be some point of $\Sigma$
whose stabilizer contains the flow $\phi^t$ of the vector field $Z$ given in Lemma~\ref{thm:Zparabolic}.
This is clearly a contradiction, because such $Z$ has no singularity.

In order to prove the claim, consider the set $\Sigma^2=\big\{x\in M:\mathrm{dim}(S_x)=2\big\}$.
This is a closed subset of $M$, because $2$ is the highest possible dimension of the stabilizers of
the $\mathds T^3$-action. If $\Sigma^2$ is non empty, then there exists an $f$-invariant measure on
$\Sigma^2$, and by Poincar\'e recurrence theorem there is at least one recurrent point $x_0\in\Sigma^2$.
The Lie algebra $\mathfrak s_{x_0}$ of $S_{x_0}$ is then $\rho(f)$-recurrent, and since
$\rho(f)$ is parabolic, by Lemma~\ref{thm:prelimlineardynamics1} (applied to the $\rho(f)$-action on
the Grassmannian of $2$-planes in $\mathcal T$), then $\mathfrak s_{x_0}$
is fixed by $\rho(f)$. There is only one $2$-plane fixed by $\rho(f)$ in $\mathcal T$ (the one spanned
by the first two vector of the basis that puts $\rho(f)$ in the normal form), and such plane
contains $Z$.

Similarly, if $\Sigma^2$ is empty, then $\Sigma^1=\big\{x\in M:\mathrm{dim}(S_x)=1\big\}$ is
closed in $M$. As above, there must be a recurrent point $x_0$ in $\Sigma^1$, and $\mathfrak s_{x_0}$
is fixed by $\rho(f)$. This implies that $\mathfrak s_{x_0}$ contains $Z$.

The proof is concluded.
\end{proof}

\begin{lem}\label{thm:everywheretimelike}
The $\T^3$-action is everywhere timelike.
\end{lem}
\begin{proof}
If not, there exists $x\in M$ such that the restriction $\mathbf g_x$ of the metric $\mathbf g$ to
$T_x\T^3x\cong\mathcal T$ is lightlike, i.e., positive semi-definite (note that the $Z$ is a lightlike
vector of such restriction, that cannot be positive definite).  Consider the $f$-invariant
compact subset $M_+=\big\{x\in M:\mathbf g_x\ \text{is positive semi-definite}\big\}$;
it has an $f$-invariant measure, and by Poincar\'e recurrence theorem there is a recurrent
point $x_0\in M_+$ for $f$. Also the metric $\mathbf g_{x_0}$ on $\mathcal T$ is $\rho(f)$-recurrent,
and by Lemma~\ref{thm:prelimlineardynamics1} (applied to the $\rho(f)$-action on
the space of quadratic forms on $\mathcal T$), $\mathbf g_{x_0}$ is fixed by $\rho(f)$.
But there exists no non zero $\rho(f)$-invariant quadratic form on $\mathcal T$ whose kernel is $Z$. This is
proved with an elementary computation using the normal form \eqref{eq:normalformrhof3space} of $\rho(f)$.
\end{proof}

\section{A general covering lemma}
\label{sec:generalcovering}
Let us now go back to the general case where $\rho(f)$ is either parabolic or hyperbolic,
and proceed with the study of the geometrical structure of $M$.
The product structure of (a finite covering of) $M$ willbe established using a general
covering result.
 \begin{prop} \label{product}
 Let $M$ be a compact manifold, and let $X$ be a non singular vector field on $M$
 of the manifold $\R\times N_0$ generating an equicontinuous flow $\phi^t$ (i.e.,
 $\phi^t$ preserves some Riemannian metric). Assume there exists a codimension $1$ foliation $\cN$
such that:
\begin{itemize}
\item[-] $\cN$ is everywhere transverse to $X$;
\item[-] $\cN$ is preserved by $\phi^t$.
\end{itemize}
Then, $X$ and $\cN$ defines a global product structure in the universal cover $\tilde{M}$. More precisely, let
$x_0\in M$ and let $N_0$ be its $\cN$-leaf. Then, the map $p:\R\times N_0\to M$
defined by $p(t,x)=\phi^tx$ is a covering.
\smallskip

A generalization is available for some group actions. Namely, consider an action of a compact Lie group
$K$ on a compact manifold $M$ such that:
\begin{itemize}
\item[-] the action is locally free (in particular, all orbits have the same dimension);
\item[-] $K$ preserves a foliation $\mathcal N$ transverse to its orbits (with a complementary dimension).
\end{itemize}
Then, for all $x_0\in M$, denoting by  by $N_0$ the
leaf of $\mathcal N$ through $x_0$, the map $p:K\times N_0\to M$ defined by $p(g,x)=gx$ is an
equivariant\footnote{The action of $K$ on $K\times N_0$ is the left multiplication on the first factor.} covering.
\end{prop}
\begin{proof}
In order to prove the first statement, consider the set of Riemannian metrics
for which $X$ and $\cN$ are orthogonal, and $X$ has norm equal to $1$. The equicontinuity assumption
implies that the $\phi_t$'s generate a precompact subgroup $\Phi$ of the diffeomorphisms group of $M$.
By averaging over the compact group $\overline\Phi$, one obtains a metric $\mathbf g_*$  on $M$
which is preserved by $\phi^t$. Now, endow $\R\times N_0$ with the product
metric, where $\R$ is endowed with the Euclidean metric $\mathrm dt^2$ and $N_0$
has the induced metric from $\mathbf g_*$. Observe that the induced metric on $N_0$ is
complete (leaves of foliations in compact manifolds have \emph{bounded geometry}, i.e.,
they are complete, have bounded curvature and injectivity radius bounded from below).
One then  observes that $p$ is a local isometry; namely, $\R\times N_0$ is complete, and therefore $p$ is a covering.

For the second statement, one can choose a left-invariant Riemannian metric $\mathbf h$ on $K$, and taking
an average on $K$ one obtains a $K$-invariant Riemannian metric $\mathbf g_*$ on $M$ such that:
\begin{itemize}
\item[(a)] the $K$-orbits and the leaves of $\mathcal N$ are everywhere orthogonal;
\item[(b)] the map $K\ni k\mapsto kx_0\in Kx_0$ is a local isometry when the orbit $Kx_0$ is endowed with the
Riemannian metric induced by $\mathbf g_*$.
\end{itemize}
As above, with such a choice the equivariant map $p:K\times N_0\to M$ defined by $p(k,x)=kx$ is a local isometry, and
since $K\times N_0$ is  complete, $p$ is a covering map.
\end{proof}

 \section{On the product structure: the hyperbolic case}
 \label{sec:hyperbolicprodstructure}
 Let us now assume that $\rho(f)$ is hyperbolic; in this section we will denote by $\T$
 the torus $\T^d\subset G_0$ given in Theorem~\ref{almost.cyclic}.
 \begin{lem}
 The orthogonal distribution  $\mathcal N$ to the $\T$-foliation is integrable.
 \end{lem}
\begin{proof}
 Let $N$ be the quotient of $M$ by the $\T$-action, and $\pi: M \to N$ the projection.
 It is a compact Riemanian orbifold. The $f$-action induces an isometry $g$ of $N$.
Consider the \emph{Levi form} (i.e., the integrability tensor of the distribution $\mathcal N$)
$l: \mathcal N \times \mathcal N \to \mathcal N^\perp$. Observe that $\mathcal N^\perp$ is the
tangent bundle of the $\T$-foliation.

Let $X$ an $Y$ be two vector fields on $N$. Suppose they are $g$-invariant: $g_*X = X$ and $g_*Y = Y$.
Let $\bar{X}$ and $\bar{Y}$ their horizontal lifts on $M$. Then, $f_*\bar{X} = \bar{X}$
and $f_*\bar{Y} = \bar{Y}$. Hence, $l(\bar{X}, \bar{Y})$ is an $f$-invariant vector field tangent to the
$\T$-foliation. However, by definition of the minimal torus $\T$, the $\rho(f)$-action on it has no invariant
vector field. This means $l(\bar{X}, \bar{Y}) = 0$.

This proof will be finished thanks to the following fact.
 \end{proof}

\begin{prop}\label{thm:ginvariantfields} Let $g$ be an isometry of a Riemannian manifold $N$.
There is an open dense set $U$ such that for any $x \in U$,
any vector $u \in T_xU$ can be extended locally to a $g$-invariant
vector field.
\end{prop}
 \begin{proof} The closure of $\{g^n, n \in \Z \}$ in the isometry group of $N$ is a compact group with
 a torus as identity connected component.
 Apply Corollary~\ref{thm:faithfulfreeaction} to conclude that, for the associated isometric
 action on $N$, the isotropy group is trivial on an open dense set $U$.
 Given $x\in U$ and $v\in T_xN$, extend first $u$ to an arbitrary smooth vector field on the slice
through $x$ of the $S$-action, and then extend to $U$ using the $S$-action.
\end{proof}

We will now prove the compactness of the leaves of $\mathcal N$.
\begin{lem} Let $N_0$ be a leaf of $\mathcal N$. Then, $N_0$ is compact.
\end{lem}
\begin{proof}
The distribution $\mathcal N$ can be seen as a connection on the $\T$-principal bundle $M \to N$.
We have just proved that this connection is flat, i.e., $\mathcal N$ is integrable, which is
equivalent to the fact that its holonomy group is discrete. The leaves will be compact if we prove that
the holonomy group is indeed finite; for $x\in N$, we will denote by $\T_x$ the fiber at $x$ of the principal
bundle $M\to N$. Recall that if $c$ is a loop at $x\in N$, then the holonomy map
$H(c):\T_x\to\T_x$  is obtained by means of horizontal lifts of $c$. It commutes with the
$\T$-action and therefore it is a translation itself. In fact, $H(c)$ can be seen as an element
of the acting torus $\T$ (and so, it is independent on the base point $x$). We have a holonomy
map $H:\pi_1(N,x)\to\T$. In fact, since $\T$ is commutative, we have canonical identification
of holonomy maps defined on different base points. In other words $H(c)=H(c^\prime)$, whence
$c$ and $c^\prime$ are freely homotopic curves.

Up to replacing $f$ by some power, we can assume that the basic Riemannian isometry $g: N \to N$
is in the identity component of $\Iso(N)$ (since this group is compact).
Therefore,  any  loop $c$ is freely homotopic to $g(c)$, and hence $H(c)=H\big(g(c)\big)$.

Now, $f$ preserves all the structure, and thus if $\tilde{c}$ is a horizontal lift of $c$, then
$f\circ\tilde{c}$ is a horizontal lift of $g(c)$.
So, $fH(c)f^{-1}=H\big(g(c)\big)$. If $g(c)$ is freely homotopic to $c$, then $H(c)$ is a fixed point
of $\rho(f)$. But we know that that $\rho(f)$ has only finitely many fixed points (by the definition of $\T$).
Therefore, the holonomy group is finite.
 \end{proof}
\noindent
Apply now Proposition~\ref{product} to deduce
that $M$ is covered by a product $\T\times N_0\to M$.
The covering is finite because $N_0$ is compact.\smallskip

Observe that we can assume the leaf $N_0$ is $f$-invariant. Indeed, the leaf $N_0$ meets all the
fibers $\T_x$, and, say, it contains $\tilde{x}\in\T_x$. So after composing with a suitable translation
$t \in \T$, i.e., replacing $f$ by $t\circ f$, we can assume that $f(\tilde{x}) \in N_0$. This implies
that $f(N_0) = N_0$. Summarizing, all things  (the $\T$-action and $f$) can be lifted to the finite
cover $T \times N_0$.

\subsection{The metric} The Lorentz metric $\mathbf g$ is not necessarily a product of the
Riemannian metric on $N_0$ by that of $\T^2$. It is true that $\T^2$ and $N$ are $\mathbf g$-orthogonal.
Also, two leaves
 $\{t\} \times N$ and $\{t^\prime\} \times N$ are isometric, via the $\T^2$-action. However, the metric
 induced on each $\T^2 \times \{n \}$ may vary with $n$.
 Observe however that one can choose a same metric for all these toral orbits, of
 course keeping the same initial group acting isometrically.
\subsection{The non-elementary case}
Let $\Gamma$ be a discrete subgroup of $\mathrm{SO}(1, k-1)$. Its limit set $L_\Gamma$ in the sphere
(boundary at infinity of the hyperbolic space $\mathbb H^{k-1}$). The group $\Gamma$ is \emph{elementary parabolic}
if $L_\Gamma$ has cardinality equal to $1$, and \emph{elementary hyperbolic} if $L_\Gamma$ has cardinality equal to
$2$. It is known that if $\Gamma$ is not elementary, then $L_\Gamma$ is infinite, and the action of $\Gamma$
on $L_\Gamma$  is \emph{minimal}, i.e., every orbit is dense, see \cite{Thurston}.

If $\Gamma$  elementary hyperbolic, then $\Gamma$ is  virtually a cyclic group, i.e.,  up to a finite index,
it consists of powers $A^n$ of the same  hyperbolic element element $A$.
If $\Gamma$ is elementary parabolic, then it is virtually  a free abelian group of rank $d \leq k-2$, i.e.,
it has a finite index subgroup isomorphic to $\Z^d$.

One fact on non-elementary groups is that they contain hyperbolic elements.
More precisely, the set of fixed points of hyperbolic elements in $L_\Gamma$ is dense in $L_\Gamma$.

All the previous considerations in the case of a hyperbolic isometry $f$, extend to the case of a
non-elementary group. We get:
\begin{theo}\label{thm:nonelementaryparabolic}
Let $(M,\mathbf g)$ be a compact Lorentz manifold with $\Iso(M,\mathbf g)$ non compact,
but $\Iso_0(M,\mathbf g)$ compact, and let $\Gamma$ be the discrete part $\Iso(M,\mathbf g)/\Iso_0(M,\mathbf g)$.
Assume that $\Iso_0(M,\mathbf g)$ has some timelike orbit. Then, there is a torus $\T^k$ contained in
$\Iso_0(M,\mathbf g)$, invariant under the action by conjugacy of $\Gamma$, and such that the $\T^k$-action is
everywhere locally free and timelike.

 The $\Gamma$-action on $\T^k$ preserves some Lorentz metric on $\T^k$, which allows one to  identify $\Gamma$
 with a discrete subgroup of $\mathrm{SO}(1,k-1)$, as well as a subgroup of $\GL(k,\Z)$.

 If $\Gamma$ is not elementary parabolic, then, up to a finite covering, $M$ splits as a topological
 product $\T^k \times N$, where $N$ is a compact Riemannian manifold.
 One can modify the original metric $\mathbf g$ along the $\T^k$ orbits, and
 get a new metric $\mathbf g^{\mathrm{new}}$  with a larger isometry group,
 $\Iso(M,\mathbf g^{\mathrm{new}})\supset\Iso(M,\mathbf g)$,
 such that $(M, \mathbf g^{\mathrm{new}})$ is a pseudo-Riemannian direct product $\T^k \times N$.
\end{theo}
 \section{On the product structure: the parabolic case}
 \label{sec:parabolicprodstructure}
As above, we have a $\T^3$-principal fibration $M \to N$ over a Riemannian orbifold $N$,
and $\mathcal N$ is seen as a connection. Let $Z$ be the Killing field (as defined above), that
is the unique vector field which commutes with $f$.

Consider the codimension 2 bundle $\mathcal L=\mathcal N\oplus\R Z$.
 \begin{lem} $\mathcal L$ is integrable.
 \end{lem}
 \begin{proof} Let $X$ and $Y$ be two vector fields tangent to $\mathcal N$. As in the proof
 above in the hyperbolic case, we can choose $X$ and $Y$ to be $f$-invariant, see
 Proposition~\ref{thm:ginvariantfields}.
 Therefore, $l(X, Y)$ is also $f$-invariant, where $l$ is the Levi form of $\mathcal N$ (not of $\mathcal L$!). But
 $Z$
 is the unique $\rho(f)$-invariant vector. Thus, $l(X, Y)$ is tangent to $\R Z$.

 Now, consider the Lie bracket $[X,Z]$. The $\T^3$-action preserves $\mathcal N$, in particular,
 $[Z, X]$ is tangent to $\mathcal N$, for any $X$ tangent to $\mathcal N$. Consequently, $\mathcal L$ is integrable.
 \end{proof}
 As in the hyperbolic case, one can prove:
\begin{lem} The leaves of $\mathcal L$ are compact.\end{lem}
One can then apply Proposition \ref{product} to $\mathcal L$ and any $2$-dimensional torus
$\T^2$ transverse to $Z$ ($\mathcal L$ is invariant by $\T^3$ and so also by any such a $\T^2$).
One obtainss that $M$ is finitely covered by $\T^2 \times L_0$.

However, since $\T^2$ is not $f$-invariant, this product is not compatible with $f$.

As in the hyperbolic case, we can choose $L_0$ invariant by $f$.
A slightly deeper analysis shows that up to a finite cover $M$ is an amalgamated product, i.e.,
a quotient $(\T^3 \times L_0)/\mathds S^1$, see Subsection~\ref{sub:amalgameted} for details.

 The metric structure of $L_0$ is that of a lightlike manifold, that is $L_0$ is endowed with a
 positive semi-definite (degenerate) metric with a $1$-dimensional null space. Here, the null
 space corresponds to the foliation defined by an $\mathds S^1$-action (the flow of the vector field  $Z$).
The circle $\mathds S^1$ acts isometrically on the lightlike $L_0$ and the Lorentz $\T^3$. A
Lorentz metric can be defined on $(L_0 \times \T^3)/\mathds S^1$.

We have proven the following:
\begin{theo}\label{thm:elementaryparabolic} Let $f$ be an isometry of a compact Lorentz manifold $(M,\mathbf g)$ such
that
the action $\rho(f)$ on the toral component of $\Iso_0(M,\mathbf g)$ is parabolic. Then,
there is a new metric on $M$ having a larger isometry group such that $M$ is the amalgamated
product of a Lorentz torus $\T^3$, and a lightlike manifold $L_0$. Both  have an isometric $\mathds S^1$-action.
The isometry $f$ is obtained by means of a isometry $h$ of $L_0$ commuting with the $\mathds S^1$-action,
and a linear isometry on the Lorentz $\T^3$.

The same statement is valid if instead of a single parabolic $f$, we have an elementary parabolic
group $\Gamma$ of rank $d$. In this case, the torus has dimension $2 + d$.
\end{theo}

\subsection{Amalgamated products}\label{sub:amalgameted}
Given any two manifolds $X$ and $Y$ carrying free (left) actions of the circle $\mathds S^1$, then
one can consider the diagonal action of $\mathds S^1$ on the product $X\times Y$:
$g(x,y)=(gx,gy)$ for all $g\in\mathds S^1$, $x\in X$ and $y\in Y$. Let $Z$ be the quotient $(X\times Y)/\mathds S^1$
of this diagonal action.
Assume that $X$ is Lorentzian,
$Y$ is Riemannian, and the action of $\mathds S^1$ in each manifold is isometric;
one can define a natural Lorentzian structure on $Z$ as follows.
Let $A\in\mathfrak X(X)$ and $B\in\mathfrak X(Y)$ be smooth vector fields tangent to the fibers
of the $\mathds S^1$-action on $X$ and on $Y$ respectively; for $(x_0,y_0)\in X\times Y$, denote by
$\big[(x_0,y_0)\big]\in Z$ the
$\mathds S^1$-orbit $\big\{(gx_0,gy_0):g\in\mathds S^1\big\}$.
The subspace $T_{x_0}X\oplus B_{y_0}^\perp$ is complementary to the one-dimensional subspace
spanned by $(A_{x_0},B_{y_0})$ in $T_{x_0}X\oplus T_{y_0}Y$. If we denote by $\mathfrak q:X\times Y\to Z$ the
projection,
then the linear map  $\mathrm d\mathfrak q_{(x_0,y_0)}:T_{x_0}X\oplus T_{y_0}Y\to T_{[(x_0,y_0)]}Z\cong
\big(T_{x_0}X\oplus T_{y_0}Y\big)\big/\mathds R\cdot(A_{x_0},B_{y_0})$
restricts to an isomorphism:
\[\mathrm d\mathfrak q_{(x_0,y_0)}:T_{x_0}X\oplus B_{y_0}^\perp\stackrel\cong\longrightarrow T_{[(x_0,y_0)]}Z.\]
A Lorentzian metric can be defined on $Z$ by requiring that such isomorphism be isometric; in
order to see that this is well defined, we need to show that this definition is independent on the choice
of $(x_0,y_0)$ in the orbit $\big[(x_0,y_0)\big]$. For $g\in\mathds S^1$, denote by $\phi_g:X\to X$ and
$\psi_g:Y\to Y$ the isometries given by the action of $g$ on $X$ and on $Y$ respectively. By differentiating
at $(x_0,y_0)$ the commutative diagram:
\[\xymatrix{X\times Y\ar[rr]^{(\phi_g,\psi_g)}\ar[dr]_{\mathfrak q}&&X\times Y\ar[ld]^{\mathfrak q}\\&Z}\]
we get a commutative diagram:
\[\xymatrix{T_{x_0}X\oplus T_{y_0}Y\ar[rrrr]^{\big((\mathrm d\phi_g)_{x_0},(\mathrm d\psi_g)_{y_0}\big)}
\ar[drr]_{\mathrm d\mathfrak q_{(x_0,y_0)}}&&&&T_{gx_0}X\oplus T_{gy_0}Y\ar[lld]^{\mathrm d\mathfrak q_{(gx_0,gy_0)}}
\\&&T_{[(x_0,y_0)]}Z}.\]
As $\big((\mathrm d\phi_g)_{x_0},(\mathrm d\psi_g)_{y_0}\big)$ carries $T_{x_0}X\oplus B_{y_0}^\perp$
onto $T_{gx_0}X\oplus B_{gy_0}^\perp$, we get the following commutative diagram of isomorphisms:
\[\xymatrix{T_{x_0}X\oplus B_{y_0}^\perp\ar[rrrr]_\cong^{\big((\mathrm d\phi_g)_{x_0},(\mathrm d\psi_g)_{y_0}\big)}
\ar[drr]^\cong_{\mathrm d\mathfrak q_{(x_0,y_0)}}&&&&T_{gx_0}X\oplus B_{gy_0}^\perp\ar[lld]_\cong^{\mathrm d\mathfrak
q_{(gx_0,gy_0)}}
\\&&T_{[(x_0,y_0)]}Z}.\]
Since $\big((\mathrm d\phi_g)_{x_0},(\mathrm d\psi_g)_{y_0}\big)$ is an isometry, the above diagram shows
that the metric induced by $\mathrm d\mathfrak q_{(x_0,y_0)}$ coincides with the metric induced by
$\mathrm d\mathfrak q_{(gx_0,gy_0)}$. This shows that the Lorentzian metric tensor on $Z$ is well defined.

As to the topology of $Z$, we have the following:
\begin{lem}\label{thm:pi1Z}
If $X$ and $Y$ are simply connected, then $Z$ is simply connected. If the product of the fundamental
groups $\pi_1(X)\times\pi_1(Y)$ is not a cyclic group, then
$Z$ is \emph{not} simply connected.
\end{lem}
\begin{proof}
The diagonal $\mathds S^1$-action on $X\times Y$ is free (and proper), and therefore the quotient
map $\mathfrak q:X\times Y\to Z$ is a smooth fibration. The thesis follows from an immediate
analysis of the long exact homotopy sequence of the fibration, that reads:
\[\mathds Z\cong\pi_1(\mathds S^1)\longrightarrow\pi_1(X)\times\pi_1(Y)\longrightarrow\pi_1(Z)
\longrightarrow\pi_0(\mathds S^1)\cong\{1\}.\qedhere\]
\end{proof}
\end{section}
\begin{section}{Proof of Corollary~\ref{thm:somewhere-everywhere} and Theorem~\ref{thm:simplyconnected}}
\begin{proof}[Proof of Corollary~\ref{thm:somewhere-everywhere}]
This is one of the steps of the proof of our structure result, see Lemma~\ref{thm:everywheretimelike}.
\end{proof}
\begin{proof}[Proof of Theorem~\ref{thm:simplyconnected}]
By the structure result of \cite{Ghani1},
compact Lorentzian manifolds admitting an isometric
action of (some covering of) $\mathrm{SL}(2,\R)$ or of an oscillator group are not simply connected. Thus,
if $M$ is simply connected, by Theorem~\ref{theorem1} $\mathrm{Iso}_0(M,\mathbf g)$ is compact.
Now, if $\mathrm{Iso}(M,\mathbf g)$ has infinitely many connected components, then (a finite covering of) $M$
is not simply connected. When $\Gamma=\mathrm{Iso}(M,\mathbf g)/\mathrm{Iso}_0(M,\mathbf g)$
is not elementary parabolic, this follows directly from Theorem~\ref{thm:nonelementaryparabolic}.
When $\Gamma$ is elementary parabolic, this follows from Theorem~\ref{thm:elementaryparabolic}
and the second statement of Lemma~\ref{thm:pi1Z}.
Hence, $\mathrm{Iso}(M,\mathbf g)$ is compact.
\end{proof}

\end{section}

\end{document}